\documentclass[12pt]{amsart}
\usepackage{amsmath, amssymb, latexsym, amsthm}
\usepackage[all]{xy}

      \newenvironment{changemargin}[2]{\begin{list}{}{
         \setlength{\topsep}{0pt}\setlength{\leftmargin}{0pt}
         \setlength{\rightmargin}{0pt}
         \setlength{\listparindent}{\parindent}
         \setlength{\itemindent}{\parindent}
         \setlength{\parsep}{0pt plus 1pt}
         \addtolength{\leftmargin}{#1}\addtolength{\rightmargin}{#2}
         }\item }{\end{list}}

\newcommand{\arx}[1]{\texttt{http://arxiv.org/abs/#1}}
\newcommand{\emp}{excluded middle property}
\newcommand{\fU}{\mathfrak{U}}
\newcommand{\fV}{\mathfrak{V}}
\newcommand{\fW}{\mathfrak{W}}
\newcommand{\sr}[2]{{\txt{$#1$\\$#2$}}}
\newcommand{\seq}[1]{\{#1\}_{n\in\N}}

\newcommand{\LE}{\preccurlyeq}
\newcommand{\GE}{\succcurlyeq}

\newcommand{\Cal}{\mathcal}

\newcommand{\A}{\forall}
\newcommand{\cA}{\mathcal{A}}
\newcommand{\B}{{\Cal B}}

\newcommand{\BG}{\B_\Gamma}
\newcommand{\BT}{\B_\Tau}
\newcommand{\BTstar}{\B_{\Tau^*}}
\newcommand{\BO}{\B_\Omega}

\newcommand{\sT}{\mathsf{T}}
\newcommand{\wX}{\mathsf{wX}}
\newcommand{\sX}{\mathsf{X}}

\newcommand{\sB}{\mathsf{B}}
\newcommand{\Tau}{\mathrm{T}}
\newcommand{\cF}{{\Cal F}}
\newcommand{\J}{{\Cal J}}

\newcommand{\M}{{\Cal M}}
\newcommand{\N}{\naturals}
\newcommand{\NN}{{{}^{\naturals}\naturals}}
\newcommand{\NZ}{{{}^{\naturals}\Z}}
\renewcommand{\inf}{P_\oo(\N)}
\renewcommand{\O}{\Cal O}

\newcommand{\R}{\reals}

\newcommand{\cU}{{\Cal U}}
\newcommand{\Union}{\bigcup}
\newcommand{\cV}{{\Cal V}}
\newcommand{\cW}{{\Cal W}}
\newcommand{\Z}{\mathbb{Z}}
\newcommand{\Impl}{\Rightarrow}
\long\def\forget#1\forgotten{}

\renewcommand{\b}{\mathfrak{b}}
\renewcommand{\t}{\mathfrak{t}}
\newcommand{\kwt}{\kappa_{\w\tau}}
\renewcommand{\c}{\mathfrak{c}}
\renewcommand{\d}{\mathfrak{d}}

\renewcommand{\i}{\item}
\renewcommand{\k}{\kappa}
\newcommand{\oo}{\infty}
\newcommand{\p}{\mathfrak{p}}
\newcommand{\fx}{\mathfrak{x}}
\newcommand{\s}{\mathfrak{s}}
\newcommand{\w}{\omega}
\newcommand{\x}{\times}

\newcommand{\Iff}{\Leftrightarrow}

\newcommand{\nin}{\not\in}

\newcommand{\sbst}{\subseteq}
\newcommand{\spst}{\supseteq}
\newcommand{\sm}{\setminus}

\newcommand{\as}{\subseteq^*}
\renewcommand{\pi}{pseudo-intersection}
\renewcommand{\(}{\left(}
\renewcommand{\)}{\right)}

\renewcommand{\>}{\rangle}
\newcommand{\E}{\exists}

\newcommand{\cov}{\mathsf{cov}}
\newcommand{\add}{\mathsf{add}}
\newcommand{\cof}{\mathsf{cof}}

\newcommand{\non}{\mathsf{non}}
\newcommand{\unif}{\mathsf{non}}

\newcommand{\impl}{\to}

\newtheorem{thm}{Theorem}[section]
\newtheorem{prop}[thm]{Proposition}
\newtheorem{prob}[thm]{Problem}
\newtheorem{lem}[thm]{Lemma}
\newtheorem{cor}[thm]{Corollary}

\theoremstyle{definition}
\newtheorem{definition}[thm]{Definition}
\theoremstyle{remark}
\newtheorem{rem}[thm]{Remark}
\newcommand{\be}{\begin{enumerate}}
\newcommand{\ee}{\end{enumerate}}
\newcommand{\bi}{\begin{itemize}}
\newcommand{\ei}{\end{itemize}}
\newcommand{\Subsection}[1]{\goodbreak\subsection{#1}}

\newcommand{\sone}{\mathsf{S}_1}    \newcommand{\sfin}{\mathsf{S}_{fin}}
\newcommand{\ufin}{\mathsf{U}_{fin}}

\newcommand{\naturals}{{\mathbb N}}
\newcommand{\reals}{{\mathbb R}}

\author{Boaz Tsaban}
\thanks{This paper constitutes a part of the author's doctoral dissertation at
Bar-Ilan University.}
\address{Department of Mathematics and Computer Science, Bar-Ilan University,
Ramat-Gan 52900, Israel}
\email{tsaban@macs.biu.ac.il, http://www.cs.biu.ac.il/\~{}tsaban}

\title{Selection principles and the Minimal Tower problem}

\begin{document}
\begin{abstract}
We study diagonalizations of covers using various selection
principles, where the covers are related to linear
quasiorderings ($\tau$-covers).
This includes: equivalences and nonequivalences,
combinatorial characterizations, critical cardinalities and
constructions of special sets of reals.
This study leads to a solution of a topological problem which was
suggested to the author by Scheepers (and stated in \cite{tau})
and is related to the Minimal Tower problem.

We also introduce a variant of the notion of $\tau$-cover,
called $\tau^*$-cover, and settle some problems for this
variant which are still open in the case of $\tau$-covers.
This new variant introduces new (and tighter) topological
and combinatorial lower bounds on the Minimal Tower problem.
\end{abstract}

\keywords{%
Gerlits-Nagy property $\gamma$-sets,
$\gamma$-cover,
$\omega$-cover,
$\tau$-cover,
tower,
selection principles,
Borel covers,
open covers}
\subjclass{03E05, 54D20, 54D80}

\maketitle

\section{Introduction}

\Subsection{Combinatorial spaces}
We consider zero-dimensional sets of real numbers.
For convenience, we may consider other spaces with more evident
combinatorial structure,
such as the \emph{Baire space} $\NN$ of infinite sequences
of natural numbers, and
the \emph{Cantor space} ${}^\N\{0,1\}$ of infinite sequences of
``bits'' (both equipped with the product topology).
The Cantor space can be identified with $P(\N)$ using
characteristic functions. We will often work in the
subspace $\inf$ of $P(\N)$, consisting of the infinite sets of natural
numbers. These spaces, as well as any separable,
zero-dimensional metric space,
are homeomorphic to sets of reals, thus our results about sets of
reals can be thought of as
talking about this more general case.

\Subsection{Selection principles}
Let $\fU$ and $\fV$ be collections of covers of a space $X$.
The following selection hypotheses have a
long history for the case when the collections
$\fU$ and $\fV$ are topologically significant.
\begin{itemize}
\item[$\sone(\fU,\fV)$:]
For each sequence $\seq{\cU_n}$ of members of $\fU$,
there is a sequence
$\seq{V_n}$ such that for each $n$ $V_n\in\cU_n$, and $\seq{V_n}\in\fV$.
\item[$\sfin(\fU,\fV)$:]
For each sequence $\seq{\cU_n}$
of members of $\fU$, there is a sequence
$\seq{\cF_n}$ such that each $\cF_n$ is a finite
(possibly empty) subset of $\cU_n$, and
$\Union_{n\in\N}\cF_n\in\fV$.
\item[$\ufin(\fU,\fV)$:]
For each sequence
$\seq{\cU_n}$ of members of $\fU$
\emph{which do not contain a finite subcover},
there exists a sequence $\seq{\cF_n}$
such that for each $n$ $\cF_n$ is a finite (possibly empty) subset of
$\cU_n$, and
$\seq{\cup\cF_n}\in\fV$.
\end{itemize}

\begin{figure}[!h]
{\tiny
$$\xymatrix@C=-2pt@R=10pt{
&
&
& \sr{\ufin(\Gamma,\Gamma)}{\b}\ar[rr]\ar@{.>}[dr]^?
&
& \sr{\ufin(\Gamma,\Omega)}{\d}\ar[rrrrr]\ar@/_/@{.>}[dl]_?
&
&
&
&
&
&
& \sr{\ufin(\Gamma,\O)}{\d}
\\
&
&
&
& \sr{\sfin(\Gamma,\Omega)}{\d}\ar[ur]
\\
& \sr{\sone(\Gamma,\Gamma)}{\b}\ar[rr]\ar[uurr]
&
& \sr{\sone(\Gamma,\Omega)}{\d}\ar[rrr]\ar[ur]
&
&
& \sr{\sone(\Gamma,\O)}{\d}\ar[uurrrrrr]
\\
  \sr{\sone(\BG,\BG)}{\b}\ar[ur]\ar[rr]
&
& \sr{\sone(\BG,\BO)}{\d}\ar[ur]\ar[rrr]
&
&
& \sr{\sone(\BG,\B)}{\d}\ar[ur]
\\
&
&
&
& \sr{\sfin(\Omega,\Omega)}{\d}\ar'[u]'[uu][uuu]
\\
\\
&
& \sr{\sfin(\BO,\BO)}{\d}\ar[uuu]\ar[uurr]
\\
& \sr{\sone(\Omega,\Gamma)}{\p}\ar'[r][rr]\ar'[uuuu][uuuuu]
&
& \sr{\sone(\Omega,\Omega)}{\cov(\M)}\ar'[uuuu][uuuuu]\ar'[rr][rrr]\ar[uuur]
&
&
& \sr{\sone(\O,\O)}{\cov(\M)}\ar[uuuuu]
\\
  \sr{\sone(\BO,\BG)}{\p}\ar[uuuuu]\ar[rr]\ar[ur]
&
& \sr{\sone(\BO,\BO)}{\cov(\M)}\ar[uu]\ar[ur]\ar[rrr]
&
&
& \sr{\sone(\B,\B)}{\cov(\M)}\ar[uuuuu]\ar[ur]
}$$
}
\caption{The surviving classes}\label{survive}
\end{figure}
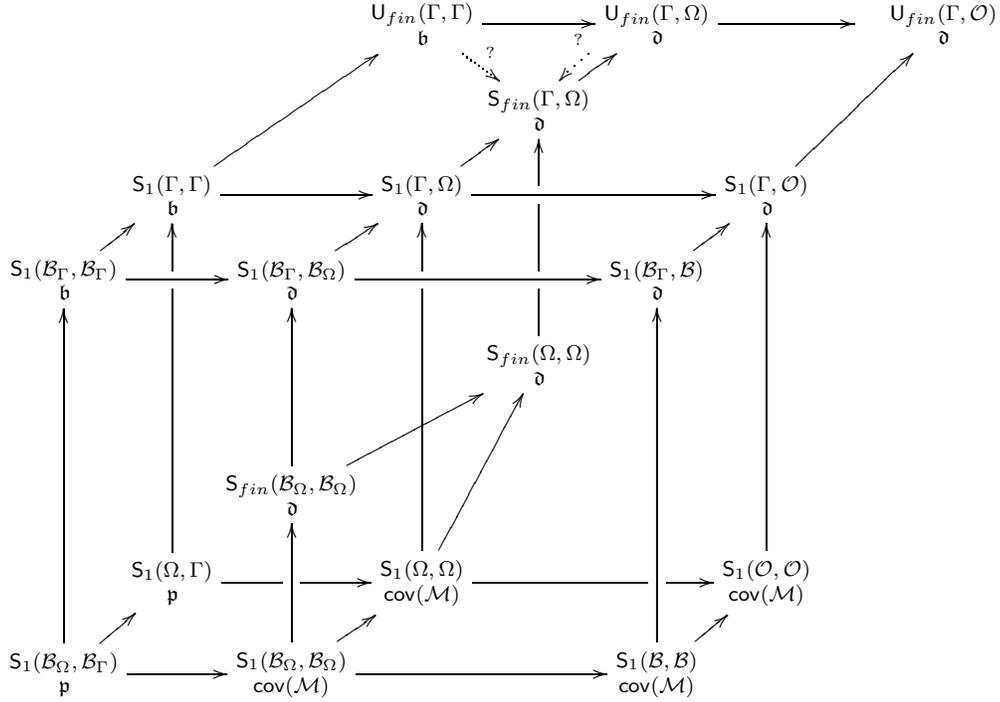

We make the convention that
\begin{quote}
\emph{The space $X$ is infinite and all covers we consider
are assumed not to have $X$ as an element.}
\end{quote}
An $\w$-cover of $X$ is a cover such that each finite subset of
$X$ is contained in some member of the cover. It is a
$\gamma$-cover if it is infinite, and each element of $X$ belongs
to all but finitely many members of the cover. Following
\cite{coc2} and \cite{CBC}, we consider the following types of
covers: \bi \i[$\O$] (respectively, $\B$): The collection of
countable open (respectively, Borel) covers of $X$. \i[$\Omega$]
(respectively, $\BO$): The collection of countable open
(respectively, Borel) $\w$-covers of $X$. \i[$\Gamma$]
(respectively, $\BG$): The collection of countable open
(respectively, Borel) $\gamma$-covers of $X$. \ei The inclusions
among these classes can be summarized as follows:
$$\begin{matrix}
\BG      & \impl & \BO      & \impl & \B \cr
\uparrow &       & \uparrow &       & \uparrow \cr
\Gamma   & \impl & \Omega   & \impl & \O
\end{matrix}$$
These inclusions and the
properties of the selection hypotheses lead to a
complicated diagram depicting how
the classes defined this way interrelate. However, only a few of
these classes are really distinct. Figure \ref{survive}
contains the distinct ones among these classes,
together with their critical
cardinalities, which were derived in
\cite{coc2} and in \cite{CBC};
see definition in Section \ref{borelimages}.
The only unsettled implications in this diagram are marked with
dotted arrows.

\Subsection{$\tau$-covers}
A cover of a space $X$ is \emph{large} if each element of
$X$ is covered by infinitely many members of the cover.
Following \cite{tau}, we consider the following type of cover.
A \emph{large} cover $\cU$ of $X$ is a
\emph{$\tau$-cover}
of $X$ if 
for each $x,y\in X$ we have either $x\in U$ implies $y\in U$ for all but finitely
many members
$U$ of the cover $\cU$, or $y\in U$ implies $x\in U$ for all but finitely many
$U\in \cU$.

A \emph{quasiordering} $\LE$ on a set $X$ is a reflexive and transitive
relation on $X$. It is \emph{linear} if for all
$x,y\in X$ we have $x\LE y$ or $y\LE x$.
A $\tau$-cover $\cU$ of a space $X$
induces a linear quasiordering
$\LE$ on $X$ by:
$$x\LE y \Iff x\in U\impl y\in U\mbox{ for all but finitely many }U\in\cU.$$
If a countable $\tau$-cover is Borel, then the induced
$\LE=\{\<x,y\> : x\LE y\}$
is a Borel subset of $X\times X$.
We let $\Tau$ and $\BT$ denote the collections of countable open and Borel
$\tau$-covers of $X$, respectively.
We have the following implications.
$$\begin{matrix}
\BG      & \impl & \BT      & \impl & \BO      & \impl & \B \cr
\uparrow &       & \uparrow &       & \uparrow &       & \uparrow\cr
\Gamma   & \impl & \Tau      & \impl & \Omega   & \impl & \O
\end{matrix}$$

There is a simple hierarchy between the selection principles:
For each $\fU,\fV$ in $\{\O,\Omega,\Tau,\Gamma\}$ or in
$\{\B,\BO,\BT,\BG\}$,
we have that $\sone(\fU,\fV)\impl\sfin(\fU,\fV)\impl\ufin(\fU,\fV)$.
The implication $\sfin(\fU,\fV)\impl\ufin(\fU,\fV)$ needs a little
care when
$\fV$ is $\Tau$ or $\BT$: It holds due to the following lemma.
\begin{lem}\label{sfinimplufin}
Assume that $\cU=\Union_{n\in\N}\cF_n$, where each $\cF_n$ is finite,
is a $\tau$-cover of a space $X$.
Then either $\cup\cF_n = X$ for some $n$, or else
$\cV=\seq{\cup\cF_n}$ is also a $\tau$-cover of $X$.
\end{lem}
\begin{proof}
Assume that $\cup\cF_n \neq X$ for all $n$.
Then, as $\cU$ is an $\w$-cover of $X$, so is $\cV$.
In particular, $\cV$ is a large cover of $X$.
Now fix any $x,y\in X$ such that $x\in U\impl y\in U$ for all but
finitely many $U\in \cU$, and let $F=\{n : (\E U\in \cF_n)\ x\in U\mbox{ and }y\nin U \}$.
Then $F$ is finite and contains the set of $n$'s such that $x\in\cup\cF_n$ and $y\nin\cup\cF_n$.
\end{proof}

\Subsection{Equivalences}\label{equivsec}

The notion of $\tau$-covers introduces seven new pairs---namely,
$(\Tau,\O)$, $(\Tau,\Omega)$, $(\Tau,\Tau)$, $(\Tau,\Gamma)$,
$(\O,\Tau)$, $(\Omega,\Tau)$, and $(\Gamma,\Tau)$---to
which any of the selection operators $\sone$, $\sfin$, and $\ufin$ can be applied.
This makes a total of $21$ new selection hypotheses.
Fortunately, some of them are easily eliminated, using the arguments of
\cite{coc1} and \cite{coc2}. We will repeat the reasoning briefly for our case.
The details can be found in the cited references.

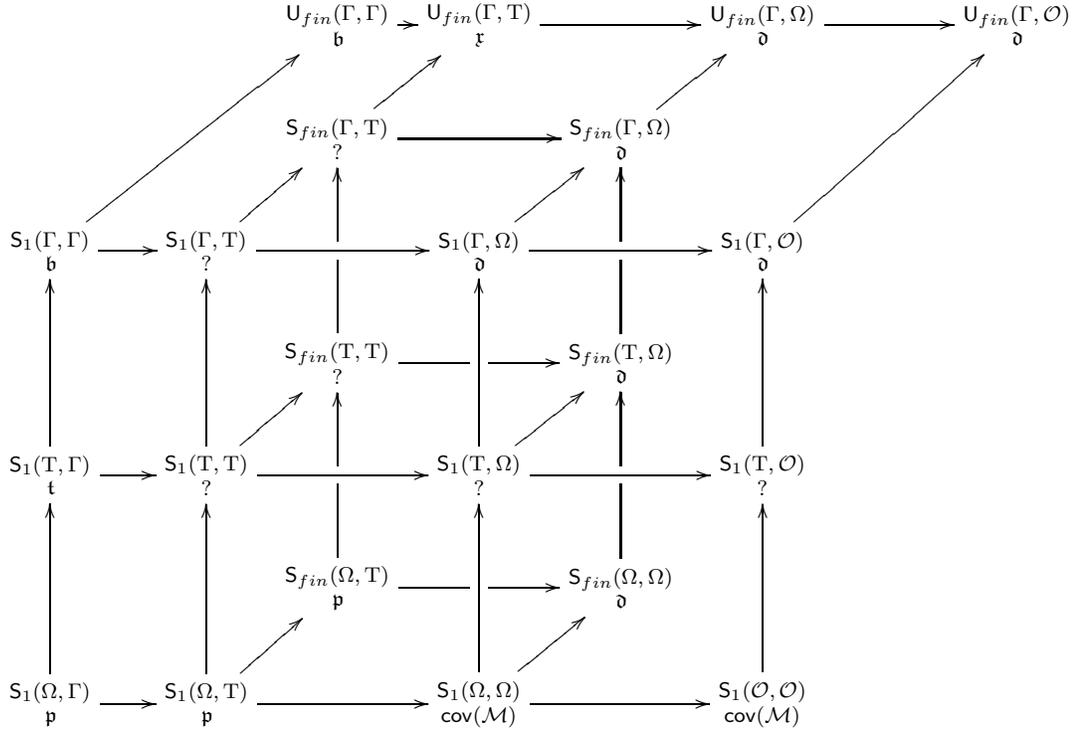
\begin{figure}[!h]
{\tiny
$$\xymatrix@C=7pt@R=20pt{
&
&
& \sr{\ufin(\Gamma,\Gamma)}{\b}\ar[r]
& \sr{\ufin(\Gamma,\Tau)}{\fx}\ar[rr]
&
& \sr{\ufin(\Gamma,\Omega)}{\d}\ar[rrrr]
&
&
&
& \sr{\ufin(\Gamma,\O)}{\d}
\\
&
&
& \sr{\sfin(\Gamma,\Tau)}{?}\ar[rr]\ar[ur]
&
& \sr{\sfin(\Gamma,\Omega)}{\d}\ar[ur]
\\
\sr{\sone(\Gamma,\Gamma)}{\b}\ar[uurrr]\ar[rr]
&
& \sr{\sone(\Gamma,\Tau)}{?}\ar[ur]\ar[rr]
&
& \sr{\sone(\Gamma,\Omega)}{\d}\ar[ur]\ar[rr]
&
& \sr{\sone(\Gamma,\O)}{\d}\ar[uurrrr]
\\
&
&
& \sr{\sfin(\Tau,\Tau)}{\mathbf{?}}\ar'[r][rr]\ar'[u][uu]
&
& \sr{\sfin(\Tau,\Omega)}{\d}\ar'[u][uu]
\\
\sr{\sone(\Tau,\Gamma)}{\t}\ar[rr]\ar[uu]
&
& \sr{\sone(\Tau,\Tau)}{\mathbf{?}}\ar[uu]\ar[ur]\ar[rr]
&
& \sr{\sone(\Tau,\Omega)}{\mathbf{?}}\ar[uu]\ar[ur]\ar[rr]
&
& \sr{\sone(\Tau,\O)}{\mathbf{?}}\ar[uu]
\\
&
&
& \sr{\sfin(\Omega,\Tau)}{\p}\ar'[u][uu]\ar'[r][rr]
&
& \sr{\sfin(\Omega,\Omega)}{\d}\ar'[u][uu]
\\
\sr{\sone(\Omega,\Gamma)}{\p}\ar[uu]\ar[rr]
&
& \sr{\sone(\Omega,\Tau)}{\p}\ar[uu]\ar[ur]\ar[rr]
&
& \sr{\sone(\Omega,\Omega)}{\cov(\M)}\ar[uu]\ar[ur]\ar[rr]
&
& \sr{\sone(\O,\O)}{\cov(\M)}\ar[uu]
}$$
}
\caption{The surviving classes for the open case}\label{survopen}
\end{figure}

First, the properties $\sone(\O,\Tau)$ and $\sfin(\O,\Tau)$ imply
$\sfin(\O,\Omega)$, and thus hold only in trivial cases (see
Section 6 of \cite{strongdiags}). Next, $\sfin(\Tau,\O)$ is
equivalent to $\ufin(\Tau,\O)$, since if the finite unions cover,
then the original sets cover as well. Now, since finite unions can
be used to turn any countable cover which does not contain a
finite subcover into a $\gamma$-cover \cite{coc2},
we have the following equivalences\footnote{%
We identify each property with the collection of sets satisfying
this property. Thus, for properties $P$ and $Q$, we may write $X\in P$,
$P\sbst Q$, etc.}:
\bi
\i $\ufin(\Tau,\Gamma)=\ufin(\Gamma,\Gamma)$,
\i $\ufin(\O,\Tau)=\ufin(\Omega,\Tau)=\ufin(\Tau,\Tau)=\ufin(\Gamma,\Tau)$,
\i $\ufin(\Tau,\Omega)=\ufin(\Gamma,\Omega)$; and
\i $\ufin(\Tau,\O)=\ufin(\Gamma,\O)$.
\ei
In Corollary \ref{tg-equiv} we get that
$\sone(\Tau,\Gamma)=\sfin(\Tau,\Gamma)$.
We are thus left with eleven new properties,
whose positions with respect to the other properties
are described in Figure \ref{survopen}.
In this Figure, as well as in the one to come,
there still exist quite many unsettled possible implications.

\goodbreak
\Subsection{Equivalences for Borel covers}
For the Borel case we have the diagram corresponding to
Figure \ref{survopen}, but
in this case, more equivalences are known \cite{CBC}:
$\sone(\BG,\BG)=\ufin(\BG,\BG)$,
$\sone(\BG,\BO)=\sfin(\BG,\BO)=\ufin(\BG,\BO)$, and
$\sone(\BG,\B)=\ufin(\BG,\B)$.
In addition, each selection principle for Borel covers
implies the corresponding selection principle for open covers.

\bigskip

This paper is divided into two parts.
Part 1 consists of Sections 2--4,
and Part 2 consists of the remaining sections.
In Section 2 we study subcover-type properties
and their applications to the study of the new
selection principles.
In Section 3 we characterize some of the properties
in terms of combinatorial properties of Borel images.
In Section 4 we find the critical cardinalities of
most of the new properties, and apply the results to
solve a topological version of the minimal
tower problem, which was suggested to us by Scheepers
and stated in \cite{tau}.

It seems that some new mathematical tools
are required to solve some of the remaining
open problems, as the special properties of
$\tau$-covers usually do not allow application
of standard methods developed during the study of
classical selection principles.
For this very reason, we believe that these are
the important problems which must be addressed
in the future.
However, we suggest in the second part of this paper
two relaxations
of the notion of $\tau$-cover, which are easier
to work with and may turn out useful in the
study of the original problems.
We demonstrate this by proving results which
are still open for the case of usual $\tau$-covers.

\part{$\tau$-covers}

\section{Subcovers with stronger properties}

\begin{definition}
Let $X$ be a set of reals, and $\fU$, $\fV$ collections of
covers of $X$. We say that $X$ satisfies $\binom{\fU}{\fV}$
(read: \emph{$\fU$ choose $\fV$}) if for each cover $\cU\in\fU$
there exists a subcover $\cV\sbst\cU$ such that $\cV\in\fV$.
\end{definition}
Observe that for any pair $\fU,\fV$ of collections of
countable covers we have that the property
$\sfin(\fU,\fV)$ implies $\binom{\fU}{\fV}$.
Gerlits and Nagy \cite{gerlitsnagy} proved that for $\fU=\Omega$
and $\fV=\Gamma$, the converse also holds, in fact,
$\sone(\Omega,\Gamma)=\binom{\Omega}{\Gamma}$. But in general
the property $\binom{\fU}{\fV}$ can be strictly weaker than
$\sfin(\fU,\fV)$.

A useful property of this notion is the following.
\begin{lem}[Cancellation Laws]\label{cancellation}
For collections of covers $\fU,\fV,\fW$,
\be
\i $\binom{\fU}{\fV}\cap\binom{\fV}{\fW}\sbst\binom{\fU}{\fW}$,
\i $\binom{\fU}{\fV}\cap\sfin(\fV,\fW)\sbst\sfin(\fU,\fW)$,
\i $\sfin(\fU,\fV)\cap\binom{\fV}{\fW}\sbst\sfin(\fU,\fW)$; and
\i $\binom{\fU}{\fV}\cap\sone(\fV,\fW)\sbst\sone(\fU,\fW)$,
\i If $\fW$ is closed under taking supersets, then
$\sone(\fU,\fV)\cap\binom{\fV}{\fW}\sbst\sone(\fU,\fW)$.
\ee
Moreover, if $\fU\spst\fV\spst\fW$, then equality holds in (1)--(5).
\end{lem}
\begin{proof}
(1) is immediate. To prove (2), we can
apply $\sfin(\fV,\fW)$ to $\fV$-subcovers of the given covers.
(4) is similar to (2).

(3) Assume that $\cU_n\in\fU$,
$n\in\N$, are given. Apply $\sfin(\fU,\fV)$ to choose finite subsets
$\cF_n\sbst\cU_n$, $n\in\N$, such that $\cV=\Union_{n\in\N}\cF_n\in\fV$.
By $\binom{\fV}{\fW}$, there exists a subset $\cW$ of $\cV$ such that
$\cW\in\fW$. Then for each $n$ $\cW\cap\cF_n$ is a finite (possibly empty)
subset of $\cU_n$, and $\Union_{n\in\N}(\cW\cap\cF_n)=\cW\in\fW$.
To prove (5), observe that the resulting cover $\cV$ contains an element
of $\fW$, and as $\fW$ is closed under taking supersets, $\cV\in\fW$
as well.

It is clear that reverse inclusion (and therefore equality)
hold in (1)--(5) when $\fU\spst\fV\spst\fW$.
\end{proof}

\begin{cor}\label{equivalences2}
Assume that $\fU\spst\fV$.
Then the following equivalences hold:
\be
\i $\sfin(\fU,\fV) = \binom{\fU}{\fV}\cap\sfin(\fV,\fU)$.
\i If $\fV$ is closed under taking supersets, then
$\sone(\fU,\fV) = \binom{\fU}{\fV}\cap\sone(\fV,\fU)$.
\ee
\end{cor}
\begin{proof}
We prove (1).
Clearly $\sfin(\fU,\fV)$ implies $\binom{\fU}{\fV}$ and $\sfin(\fV,\fU)$.
On the other hand, by applying the Cancellation Laws (2) and then (3) we have that
$$\binom{\fU}{\fV}\cap\sfin(\fV,\fU)\sbst\sfin(\fU,\fU)\cap\binom{\fU}{\fV}\sbst\sfin(\fU,\fV).$$
\end{proof}

\Subsection{When every $\tau$-cover contains a $\gamma$-cover}


\begin{thm}\label{equivalences}
The following equivalences hold:
\be
\i $\sone(\Tau,\Gamma)=\binom{\Tau}{\Gamma}\cap\sfin(\Gamma,\Tau)$,
\i $\sone(\BT,\BG)=\binom{\BT}{\BG}\cap\sfin(\BG,\BT)$,
\ee
\end{thm}
\begin{proof}
(1) By the Cancellation Laws \ref{cancellation},
$\binom{\Tau}{\Gamma}\cap\sfin(\Gamma,\Tau)\sbst
\sfin(\Gamma,\Gamma)$.
In \cite{coc2} it was proved that $\sfin(\Gamma,\Gamma)=\sone(\Gamma,\Gamma)$.
Thus, $\binom{\Tau}{\Gamma}\cap\sfin(\Gamma,\Tau)\sbst
\binom{\Tau}{\Gamma}\cap\sone(\Gamma,\Gamma)$, which by the Cancellation Laws
is a subset of $\sone(\Tau,\Gamma)$. The other direction is immediate.
\forget
(1) We need to show that
$\binom{\Tau}{\Gamma}\cap\sfin(\Gamma,\Tau)\sbst\sone(\Tau,\Gamma)$.
Assume that $X\in\binom{\Tau}{\Gamma}\cap\sfin(\Gamma,\Tau)$, and
that $\cU_n$, $n\in\N$, are countable open $\tau$-covers of $X$.
Using the property $\binom{\Tau}{\Gamma}$, we may assume that the
covers $\cU_n$ are $\gamma$-covers of $X$. Thus it suffices to
show that $X$ satisfies $\sone(\Gamma,\Gamma)$. Using the
$\sfin(\Gamma,\Tau)$ property, we can choose finite subsets
$\cF_n\sbst\cU_n$ such that $\cU=\Union_{n\in\N}\cF_n$ is a
$\tau$-cover of $X$. Applying the property $\binom{\Tau}{\Gamma}$
again, we have that $\cU$ contains a $\gamma$-cover of $X$. This
implies that $X$ satisfies $\sfin(\Gamma,\Gamma)$. But
$\sfin(\Gamma,\Gamma)=\sone(\Gamma,\Gamma)$ \cite{coc2}.
%
\forgotten

(2) is similar.
\end{proof}

\begin{cor}\label{tg-equiv}
The following equivalences hold:
\be
\i $\sone(\Tau,\Gamma)=\sfin(\Tau,\Gamma)$;
\i $\sone(\BT,\BG)=\sfin(\BT,\BG)$.
\ee
\end{cor}

Using similar arguments, we have the following.
\begin{thm}\label{equiv2}
The following equivalences hold:
\be
\i $\sone(\Omega,\Gamma) = \binom{\Tau}{\Gamma}\cap\sfin(\Omega,\Tau)$;
\i $\sone(\BO,\BG) = \binom{\BT}{\BG}\cap\sfin(\BO,\BT)$.
\ee
\end{thm}

\Subsection{When every $\w$-cover contains a $\tau$-cover}

\begin{thm}\label{great}
The following inclusions hold:
\be
\i $\binom{\Omega}{\Tau}\sbst\sfin(\Gamma,\Tau)$.
\i $\binom{\BO}{\BT}\sbst\sfin(\BG,\BT)$.
\ee
\end{thm}
\begin{proof}
We will prove (1) (the proof of (2) is identical). Assume that $X$
satisfies $\binom{\Omega}{\Tau}$. If $X$ is countable then it
satisfies all of the properties mentioned in this paper. Otherwise
let $x_n$, $n\in\N$, be distinct elements in $X$. Assume that
$\cU_n=\{U^n_m\}_{m\in\N}$, $n\in\N$, are open $\gamma$-covers of
$X$. Define $\tilde\cU_n=\{U^n_m\sm\{x_n\}\}_{m\in\N}$. Then
$\cU=\Union_{n\in\N}\tilde\cU_n$ is an open $\w$-cover of $X$, and
thus contains a $\tau$-cover $\cV$ of $X$. Let $\LE$ be the
induced quasiordering.
\begin{lem}\label{least}
If $\<X,\LE\>$ has a least element, then $\cV$ contains a
$\gamma$-cover of $X$.
\end{lem}
\begin{proof}
Write $\cV=\{V_n\}_{n\in\N}$. Let $x_0$ be a
least element in $\<X,\LE\>$. Consider the subsequence
$\{V_{n_k}\}_{k\in\N}$ consisting
of the elements $V_n$ such that $x_0\in V_n$. Since $\tau$-covers are large,
this sequence is infinite.
For all $x\in X$ we have $x_0\LE x$, thus $x_0\in V_n\impl x\in V_n$ for
all but finitely many $n$.
Since $x_0\in V_{n_k}$ for all $k$, we have that for all but
finitely many $k$, $x\in V_{n_k}$.
\end{proof}
There are two cases to consider.

\textbf{Case 1.} For some $n$ $x_n$ is a least element in
$\<X,\LE\>$. Then $\cV$ contains a $\gamma$-cover $\tilde\cV$ of
$X$. In this case, for all $n$ $x_n$ belongs to all but finitely
many members of $\tilde\cV$, thus $\tilde\cV\cap\tilde\cU_n$ is
finite for each $n$, and $\cW=\{U : (\E
n)U\sm\{x_n\}\in\tilde\cV\}$ is a $\gamma$-cover of $X$.

\textbf{Case 2.} For each $n$ there exists $x\neq x_n$ with $x\LE
x_n$. For each $n$, $\cU_n$ is a $\gamma$-cover of $X$, thus $x$
belongs to all but finitely many members of $\cV\cap\tilde\cU_n$.
Since $x_n$ does not belong to any of the members in
$\cV\cap\tilde\cU_n$, $\cV\cap\tilde\cU_n$ must be finite. Thus,
$\cW=\{U : (\E n)U\sm\{x_n\}\in\tilde\cV\}$ is a $\tau$-cover of
$X$.
\end{proof}

If $\binom{\Omega}{\Tau}\sbst\sfin(\Tau,\Omega)$, then by Corollary \ref{equivalences2}
$\binom{\Omega}{\Tau}=\sfin(\Omega,\Tau)$.
\begin{prob}
Is $\binom{\Omega}{\Tau}=\sfin(\Omega,\Tau)$?
\end{prob}

\section{Combinatorics of Borel images}\label{borelimages}

In this section we characterize several properties in terms of
Borel images in the spaces $\NN$ and $\inf$, using the combinatorial
structure of these spaces.

\Subsection{The combinatorial structures}
A quasiorder $\le^*$ is defined on the Baire space $\NN$ by eventual dominance:
$$f\le^* g \mbox{\quad if } f(n)\le g(n) \mbox{ for all but finitely many }n.$$

A subset $Y$ of $\NN$ is called \emph{unbounded} if it is unbounded with
respect to $\le^*$.
$Y$ is \emph{dominating} if it is cofinal in $\NN$ with respect to $\le^*$,
that is, for each $f\in \NN$ there exists $g\in Y$ such that
$f\le^* g$.
$\b$ is the minimal size of
an unbounded subset of $\NN$, and $\d$ is the minimal size
of a dominating subset of $\NN$.

Define a quasiorder $\as$ on $\inf$ by
$a\as b$ if $a\sm b$ is finite.
An infinite set $a\sbst\N$ is a \emph{pseudo-intersection}
of a family $Y\sbst\inf$
if for each $b\in Y$, $a\as b$.
A family $Y\sbst\inf$ is a \emph{tower} if it is linearly quasiordered
by $\as$, and it has no \pi{}.
$\t$ is the minimal size of a tower.

A family $Y\sbst\inf$ is \emph{centered} if the intersection of
each (nonempty) finite subfamily of $Y$ is infinite.
Note that every tower in $\inf$ is centered.
A centered family $Y\sbst\inf$ is a \emph{power} if it
does not have a \pi{}.
$\p$ is the minimal size of a power.

\Subsection{The property $\binom{\BT}{\BG}$}
For a set of reals $X$ and a topological space $Z$,
we say that $Y$ is a Borel image of $X$ in $Z$ if there
exists a Borel function $f:X\to Z$ such that $f[X]=Y$.
The following classes of sets were introduced in \cite{pawlikowskireclaw}:
\begin{itemize}
\item[$\mathsf{P}$:]{ The set of $X\subseteq\reals$ such that no Borel
   image of $X$ in $\inf$ is a power,}
\item[$\mathsf{B}$:]{ The set of $X\subseteq\reals$ such that every Borel
   image of $X$ in $\NN$ is bounded
   (with respect to eventual domination);}
\item[$\mathsf{D}$:]{ The set of $X\subseteq \reals$ such that no Borel
   image of $X$ in $\NN$ is dominating.}
\end{itemize}

For a collection $\J$ of separable metrizable spaces,
let $\unif(\J)$ denote the minimal cardinality of a
separable metrizable space which is not a member of $\J$.
We also call $\non(\J)$ \emph{the critical cardinality} of the class $\J$.
The critical cardinalities of the above classes are
$\p$, $\b$, and $\d$, respectively.
These classes have the interesting property that they transfer the
cardinal inequalities $\p\le\b\le\d$ to the
inclusions $\mathsf{P}\sbst\sB\sbst\mathsf{D}$.

\begin{definition}
For each countable cover of $X$ enumerated bijectively as
$\cU=\seq{U_n}$ we associate a function
$h_\cU: X\to P(\N)$, defined
by $h_\cU(x) = \{ n : x\in U_n\}$.
\end{definition}
$\cU$ is a large cover of $X$ if, and only if, $h_\cU[X]\sbst\inf$.
As we assume that $X$ is infinite and is not a member of any of our covers,
we have that each $\w$-cover of $X$ is a large cover of $X$.
The following lemma is a key observation
for the rest of this section. Note that $h_\cU$ is a Borel function
whenever $\cU$ is a Borel cover of $X$, and $h_\cU$ is continuous
whenever all elements of $\cU$ are \emph{clopen}.
\begin{lem}[\cite{tau}]\label{notions}
Assume that $\cU$ is a countable large cover of $X$. \be \i $\cU$
is an $\w$-cover of $X$ if, and only if, $h_\cU[X]$ is centered.
\i $\cU$ contains a $\gamma$-cover of $X$ if, and only if,
$h_\cU[X]$ has a \pi{}. \i $\cU$ is a $\tau$-cover of $X$ if, and
only if, $h_\cU[X]$ is linearly quasiordered by $\as$. \ee
Moreover, if $f:X\to P(\N)$ is any function, and $\cA=\seq{O_n}$
is the clopen cover of $P(\N)$ such that $x\in O_n\Iff n\in x$,
then for $\cU=\seq{f^{-1}[O_n]}$ we have that $f=h_\cU$.
\end{lem}
This Lemma implies that $\mathsf{P}=\binom{\BO}{\BG}$ \cite{CBC}.
\begin{cor}\label{P-intersect}
$\mathsf{P}=\binom{\BO}{\BT}\cap\binom{\BT}{\BG}$.
\end{cor}

It is natural to define the following notion.
\bi
\item[$\sT$:]{
The set of $X\subseteq\reals$ such that no Borel image of $X$
in $\inf$ is a tower.}
\ei

\begin{thm}\label{T-tower}
$\sT=\binom{\BT}{\BG}$.
\end{thm}
\begin{proof}
See \cite{tau} for the \emph{clopen} version of this theorem
(a straightforward usage of Lemma \ref{notions}).
The proof for the Borel case is similar.
\end{proof}

\begin{cor}\label{T-t}
$\non(\binom{\BT}{\BG})=\non(\binom{\Tau}{\Gamma})=\t$.
\end{cor}
\begin{proof}
By Theorem \ref{T-tower}, $\t\le\non(\binom{\BT}{\BG})$. In
\cite{tau} we defined $\mathcal{T}$ to be the collection of sets
for which every countable \emph{clopen} $\tau$-cover contains a
$\gamma$-cover, and showed that $\non(\mathcal{T})=\t$. But
$\binom{\BT}{\BG}\sbst\binom{\Tau}{\Gamma}\sbst\mathcal{T}$.
\end{proof}

Clearly, $\mathsf{P}\sbst \sT$.
The cardinal inequality $\p\le\t\le\b$ suggests
pushing this further by showing that $\sT\sbst \sB$;
unfortunately this is false. Sets which are continuous images
of Borel sets are called \emph{analytic}.
\begin{thm}\label{baire-T}
Every analytic set satisfies $\sT$. In particular, $\NN\in\sT$.
\end{thm}
\begin{proof}
According to \cite{tau}, no continuous image of an analytic set
is a tower. In particular, towers are not analytic subsets of $\inf$.
Since Borel images of analytic sets are again analytic sets,
we have that every analytic set satisfies $\sT$.
\end{proof}

The following equivalences hold \cite{CBC}:
\bi
\i $\sone(\BG,\BG)=\sB$,
\i $\sone(\BG,\B)=\mathsf{D}$;
\i $\sone(\BO,\BG)=\mathsf{P}$.
\ei
Theorem \ref{baire-T} rules out an identification of $\sT$
with any of the selection principles.
However, we get the following characterization of
$\sone(\BT,\BG)$ in terms of Borel images.
\begin{thm}
$\sone(\BT,\BG) = \sT\cap\sB$.
\end{thm}
\begin{proof}
By the Cancellation Laws and Theorem \ref{T-tower},
$\sone(\BT,\BG) = \binom{\BT}{\BG}\cap\sone(\BG,\BG)=\sT\cap\sB$.
\end{proof}

\Subsection{The property $\binom{\BO}{\BT}$}
For a subset $Y$ of $\inf$ and $a\in\inf$, define
$$Y\restriction a = \{y\cap a : y\in Y\}.$$
If all sets in $Y\restriction a$ are infinite, we
say that $Y\restriction a$ is a
\emph{large restriction} of $Y$.

\begin{thm}\label{w-tau}
For a set $X$ of real numbers, the following are equivalent:
\be
\i $X$ satisfies $\binom{\BO}{\BT}$
\i For each Borel image $Y$ of $X$ in $\inf$,
if $Y$ is centered, then
there exists a large restriction of $Y$
which is linearly quasiordered by $\as$.
\ee
\end{thm}
\begin{proof}
$1\Impl 2$: Assume that $\Psi:X\to\inf$ is a Borel function, and let
$Y=\Psi[X]$. Assume that $Y$ is centered, and
consider the collection $\cA=\seq{O_n}$ where
$O_n=\{a : n\in a\}\cap Y$ for each $n\in\N$.
If the set $a = \{n : Y = O_n\}$ is infinite,
then $a$ is a \pi{} of $Y$ and we are done.
Otherwise, by removing finitely many elements from $\cA$
we get that $\cA$ is an $\w$-cover of $Y$.

Setting $U_n=\Psi^{-1}[O_n]$ for each $n$,
we have that $\cU=\seq{U_n}$ is a Borel
$\w$-cover of $X$,
which thus contains a $\tau$-cover $\seq{U_{a_n}}$ of $X$.
Let $a=\seq{a_n}$, and define a cover $\cV=\seq{V_n}$ of $X$ by
$$V_n=\begin{cases}
U_n & n\in a\cr
\emptyset & $otherwise$\cr
\end{cases}$$
Then $\cV$ is a $\tau$-cover of $X$, and by
Lemma \ref{notions},
$\Psi[X]\restriction a = h_\cU[X]\restriction a = h_{\cV}[X]$ is
linearly quasiordered by $\as$.

$2\Impl 1$: Assume that $\cU=\seq{U_n}$ is an $\w$-cover of $X$.
By Lemma \ref{notions},
$h_\cU[X]$ is centered. Let $a=\seq{a_n}$ be a large restriction of
$h_\cU[X]$ which is linearly quasiordered by $\as$, and define
$\cV$ as in $1\Impl 2$. Then $h_\cU[X]\restriction a=h_\cV[X]$.
Thus all elements in $h_\cV[X]$ are infinite
(i.e., $\cV$ is a large cover of $X$), and
$h_\cV[X]$ is linearly quasiordered by $\as$ (i.e.,
$\cV$ is a $\tau$-cover of $X$).
Then $\cV\sm\{\emptyset\}\sbst\cU$ is a $\tau$-cover of $X$.
\end{proof}

\begin{rem}\label{clopen-Pwt}
Replacing the ``Borel sets'' by ``clopen sets'' and ``Borel functions''
by ``continuous functions'' in the last proof we get that the following
properties are equivalent for a set $X$ of reals:
\be
\i Every countable \emph{clopen} $\w$-cover of
$X$ contains a $\tau$-cover of $X$.
\i For each continuous image $Y$ of $X$ in $\inf$,
if $Y$ is centered, then
there exists a large restriction of $Y$
which is  linearly quasiordered by $\as$.
\ee
We do not know whether the \emph{open} version of this result is true.
\end{rem}

\Subsection{The property $\ufin(\BG,\BT)$}

\begin{definition}
A family $Y\sbst\NN$ satisfies the
\emph{excluded middle} property if there exists $g\in\NN$
such that:
\be
\i For all $f\in Y$, $g\not\le^* f$;
\i For each $f,h\in Y$, one of the
situations $f(n) < g(n) \le h(n)$
or $h(n) < g(n) \le f(n)$ is possible only
for finitely many $n$.
\ee
\end{definition}

\begin{thm}\label{excludedmiddle}
For a set $X$ of real numbers, the following are equivalent:
\begin{enumerate}
\item{$X$ satisfies $\ufin(\BG,\BT)$;}
\item{Every Borel image of $X$ in $\NN$ satisfies the excluded
middle property.}
\end{enumerate}
\end{thm}
\begin{proof}
$1\Rightarrow 2$: For each $n$, the collection ${\Cal U}_n=\{U^n_m
: m\in\N\}$, where $U^n_m = \{f\in\NN : f(n)\le m\}$, $m\in\N$, is
an open $\gamma$-cover of $\NN$. Assume that $\Psi$ is a Borel
function from $X$ to $\NN$. By standard arguments we may assume
that $\Psi^{-1}[U^n_m]\neq X$ for all $n$ and $m$. Then the
collections $\cU_n = \{\Psi^{-1}[U^n_m] : m\in\N\}$, $n\in\N$, are
Borel $\gamma$-covers of $X$. By $\ufin(\BG,\BT)$, there exist
finite sets $F_n\sbst\N$, $n\in\N$, such that $\cU=\{\Union_{m\in
F_n}\Psi^{-1}[U^n_m] : n\in\N\}$ is a $\tau$-cover of $X$. Let
$A=\{n : F_n\neq\emptyset\}$. Note that for each $n\in A$,
$\Union_{m\in F_n}\Psi^{-1}[U^n_m]=\Psi^{-1}[U^n_{\max F_n}]$.

Define $g\in\NN$ by
$$g(n) = \begin{cases}
\max F_n + 1 & n\in A\\
0 & \mbox{otherwise}
\end{cases}$$
For all $x\in X$, as $\cU$ is a large cover of $X$,
there exist infinitely many $n\in A$ such that
$\Psi(x)\in U^n_{\max F_n}$ (that is,
$\Psi(x)(n) < g(n)$).
Let $\LE$ be the linear quasiordering of $X$ induced by the
$\tau$-cover $\cU$. Then for all $x,y\in X$, either
$x\LE y$ or $y\LE x$. In the first case we get that
for all but finitely many $n$
$\Psi(x)(n) < g(n)\to \Psi(y)(n) < g(n)$,
and in the second case we get the same assertion with $x$ and
$y$ swapped. This shows that $\Psi[X]$ satisfies the excluded
middle property.

$2\Rightarrow 1$:
Assume that $\cU_n=\{U^n_m : m\in\N\}$, $n\in\N$, are Borel covers of $X$
which do not contain a finite subcover.
Replacing each $U^n_m$ with the Borel set $\Union_{k\le m} U^n_k$ we may assume
that the sets $U^n_m$ are monotonically increasing with $m$.
Define a function $\Psi$ from $X$ to $\NN$ so that for each $x$ and $n$:
$$\Psi(x)(n) = \min\{m : x\in U^n_m\}.$$
Then $\Psi$ is a Borel map, and so $\Psi[X]$
satisfies the excluded middle property.
Let $g\in\NN$ be a witness for that.
Then $\cU=\{U^n_{g(n)-1} : n\in\N,\ g(n)>0\}$ is a $\tau$-cover of $X$:
For each $x\in X$ we have that $g\not \le^*\Psi(x)$, thus
$\cU$ is a large cover of $X$.
Moreover, for all $x,y\in X$, we have by the excluded middle property that
at least one of the assertions $\Psi(x)(n)< g(n) \le \Psi(y)(n)$
or $\Psi(y)(n)< g(n)\le\Psi(y)(n)$ is possible only
for finitely many $n$.
Then the first assertion implies that
$x\LE y$, and the second implies $y\LE x$ with respect to $\cU$.
\end{proof}

\begin{rem}\label{ufinGammaTauchar}
The analogue \emph{clopen} version of Theorem \ref{excludedmiddle}
also holds.
We do not know whether there exist an
analogue characterization of $\ufin(\Gamma,\Tau)$ (the \emph{open} version)
in terms of continuous images.
\end{rem}

\section{Critical cardinalities}

\begin{thm}\label{nontw}
$\non(\sfin(\BT,\BO))=\non(\sfin(\Tau,\Omega))=\d$.
\end{thm}
\begin{proof}
$\sfin(\BO,\BO)\sbst\sfin(\BT,\BO)\sbst\sfin(\Tau,\Omega)\sbst\sfin(\Gamma,\Omega)$,
and according to \cite{coc2} and \cite{CBC},
$\non(\sfin(\BO,\BO))=\non(\sfin(\Gamma,\Omega))=\d$.
\end{proof}

\begin{thm}\label{TGt}
$\non(\sone(\BT,\BG))=\non(\sone(\Tau,\Gamma))=\t$.
\end{thm}
\begin{proof}
By Theorem \ref{equivalences},
$\sone(\Tau,\Gamma)=\binom{\Tau}{\Gamma}\cap\sone(\Gamma,\Gamma)$, thus
by Corollary \ref{T-t},
$\non(\sone(\Tau,\Gamma))=\min\{\non(\binom{\Tau}{\Gamma}),\non(\sone(\Gamma,\Gamma))\}=
\min\{\t,\b\}=\t$.
The proof for the Borel case is similar.
\end{proof}

\begin{definition}
$\fx$ is the minimal cardinality of a family $Y\sbst\NN$
which does not satisfy the excluded middle property.
\end{definition}

Therefore $\b\le\fx\le\d$.
A family $Y\sbst\inf$ is \emph{splitting} if for each infinite $a\sbst\N$
there exists $s\in Y$ which \emph{splits} $a$, that is,
such that the sets $a\cap s$ and $a\sm s$ are infinite.
$\s$ is the minimal size of a splitting family.
In \cite{ShTb} it is proved that $\fx=\max\{\s,\b\}$.

\begin{thm}\label{non=e}
$\non(\ufin(\BG,\BT))=\non(\ufin(\Gamma,\Tau))=\fx$.
\end{thm}
\begin{proof}
By Theorem \ref{excludedmiddle}, $\non(\ufin(\BG,\BT))=\fx$. Thus,
our theorem will follow from the inclusion
$\ufin(\BG,\BT)\sbst\ufin(\Gamma,\Tau)$ once we prove that
$\non(\ufin(\Gamma,\Tau))\le\fx$. To this end, consider a family
$Y\sbst\NN$ of size $\fx$ which does not satisfy the excluded
middle property, and consider the monotone $\gamma$-covers ${\Cal
U}_n$, $n\in\N$, of $\NN$ defined in the proof of Theorem
\ref{excludedmiddle}. Then, as in that proof, we cannot extract
from these covers a $\tau$-cover of $Y$. Thus, $Y$ does not
satisfy $\ufin(\Gamma,\Tau)$.
\end{proof}

\begin{definition}
Let $\kwt$ be the minimal cardinality of a centered set $Y\sbst\inf$
such that for no $a\in\inf$, the restriction
$Y\restriction a$ is large and linearly quasiordered by $\as$.
\end{definition}

It is easy to see (either from the definitions or
by consulting the involved selection properties)
that $\kwt\le\d$ and $\p=\min\{\kwt,\t\}$.
In \cite{ShTb} it is proved that in fact
$\kwt=\p$.

\begin{lem}\label{kwt}
$\non(\binom{\BO}{\BT})=\non(\binom{\Omega}{\Tau})=\p$.
\end{lem}
\begin{proof}
Let $\mathcal{P}_{\w\tau}$ denote the property that every
\emph{clopen} $\w$-cover contains a $\gamma$-cover. Then
$\binom{\BO}{\BT}\sbst\binom{\Omega}{\Tau}\sbst\mathcal{P}_{\w\tau}$.
By Theorem \ref{w-tau} and Remark \ref{clopen-Pwt},
$\non(\binom{\BO}{\BT})=\non(\mathcal{P}_{\w\tau})=\kwt=\p$.
\end{proof}


\begin{thm}\label{non=kwt}
$\non(\sfin(\BO,\BT))=\non(\sfin(\Omega,\Tau))=\p$.
\end{thm}
\begin{proof}
By Corollary \ref{equivalences2} and Theorem \ref{nontw},
\begin{eqnarray}\nonumber
\non(\sfin(\Omega,\Tau))
& = & \min\{\non(\binom{\Omega}{\Tau}),\non(\sfin(\Tau,\Omega))\} = \cr
& = & \min\{\kwt,\d\} = \kwt = \p.
\end{eqnarray}
The proof for the Borel case is the same.
\end{proof}

\begin{cor}
$\non(\sone(\BO,\BT))=\non(\sone(\Omega,\Tau))=\p$.
\end{cor}
\begin{proof}
$\sone(\BO,\BG)\sbst\sone(\BO,\BT)\sbst\sone(\Omega,\Tau)\sbst\binom{\Omega}{\Tau}$.
\end{proof}

\section{Topological variants of the Minimal Tower problem}

Let $\c$ denote the size of the continuum. The following inequalities
are well known \cite{vD}:
$$\p\le\t\le \b\le\d\le\c.$$
For each pair except $\p$ and $\t$, it is well known that a strict
inequality is consistent.
\begin{prob}[\emph{Minimal Tower}]
Is it provable that $\p=\t$?
\end{prob}
This is one of the major and oldest problems of infinitary
combinatorics. Allusions to this problem can be
found in Rothberger's works (see, e.g., \cite{ROTH2}).

We know that $\sone(\Omega,\Gamma)\sbst\sone(\Tau,\Gamma)$,
and that $\non(\sone(\Omega,\Gamma))=\p$, and $\non(\sone(\Tau,\Gamma))=\t$.
Thus, if $\p<\t$ is consistent, then it is consistent that
$\sone(\Omega,\Gamma)\neq\sone(\Tau,\Gamma)$.
Thus the following problem, which was suggested to us by Scheepers,
is a logical lower bound on the difficulty of the
Minimal Tower problem.

\begin{prob}[\cite{tau}]\label{topoprob}
Is it consistent that $\sone(\Omega,\Gamma)\neq\sone(\Tau,\Gamma)$?
\end{prob}

We also have a Borel variant of this problem.

\begin{prob}\label{topoprob2}
Is it consistent that $\sone(\BO,\BG)\neq\sone(\BT,\BG)$?
\end{prob}

We will solve both of these problems.

For a class $\J$ of sets of real numbers with $\cup\J\nin\J$,
the \emph{additivity number}
of $\J$ is the minimal cardinality of a collection $\cF\sbst\J$ such that
$\cup\cF\nin\J$. The additivity number of $\J$ is denoted $\add(\J)$.

\begin{lem}[\cite{tau}]\label{p.i.pieces}
Assume that $Y\sbst\inf$ is linearly ordered by $\as$, and
for some $\k<\t$, $Y=\Union_{\alpha<\k}Y_\alpha$ where each
$Y_\alpha$ has a \pi{}.
Then $Y$ has a \pi{}.
\end{lem}

\begin{thm}\label{addT}
$\add(\binom{\BT}{\BG})=\add(\binom{\Tau}{\Gamma})=\t$.
\end{thm}
\begin{proof}
By Theorem \ref{T-tower} and Lemma \ref{p.i.pieces}, we have that
$\add(\binom{\BT}{\BG})=\add(\sT)=\t$. The proof that
$\add(\binom{\Tau}{\Gamma})=\t$ is not as elegant and requires a
back-and-forth usage of Lemma \ref{notions}. Assume that $\k<\t$,
and let $X_\alpha$, $\alpha<\k$, be sets satisfying
$\binom{\Tau}{\Gamma}$. Let $\cU$ be a countable open $\tau$-cover
of $X=\Union_{\alpha<\k}X_\alpha$. Then
$h_\cU[X]=\Union_{\alpha<\k}h_\cU[X_\alpha]$ is linearly
quasiordered by $\as$. Since each $X_\alpha$ satisfies
$\binom{\Tau}{\Gamma}$, for each $\alpha$ $\cU$ contains a
$\gamma$-cover of $X_\alpha$, that is, $h_\cU[X_\alpha]$ has a
pseudo-intersection. By Lemma \ref{p.i.pieces}, $h_\cU[X]$ has a
pseudo-intersection, that is, $\cU$ contains a $\gamma$-cover of
$X$.
\end{proof}

\begin{thm}\label{addS1BTBG}
$\add(\sone(\BT,\BG))=\t$.
\end{thm}
\begin{proof}
By Theorem \ref{equivalences},
$\sone(\BT,\BG) = \binom{\BT}{\BG}\cap\sone(\BG,\BG)$, and
according to \cite{AddQuad}, $\add(\sone(\BG,\BG))=\b$.
By Theorem \ref{addT}, we get that
$$\add(\sone(\BT,\BG))\ge\min\{\t,\b\}=\t.$$
On the other hand, by Theorem \ref{TGt} we have
$$\add(\sone(\BT,\BG))\le\non(\sone(\BT,\BG))=\t.$$
\end{proof}

In \cite{wqn} Scheepers proves that $\sone(\Gamma,\Gamma)$ is
closed under taking unions of size less than the distributivity
number $\mathfrak{h}$. Consequently, we get that
$\add(\sone(\Tau,\Gamma))=\t$ \cite{AddQuad}.
As it is consistent that $\sone(\Omega,\Gamma)$ is not closed under
taking finite unions \cite{GM}, we get a positive solution to
Problem \ref{topoprob}. We will now prove something stronger:
Consistently, no class between $\sone(\BO,\BG)$ and $\binom{\Omega}{\Tau}$
(inclusive) is closed under taking finite unions.
This solves Problem \ref{topoprob} as well as Problem \ref{topoprob2}.
\begin{thm}[CH]\label{2-add}
There exist sets of reals $A$ and $B$ satisfying $\sone(\BO,\BG)$,
such that $A\cup B$
does not satisfy $\binom{\Omega}{\Tau}$.
In particular, $\sone(\Tau,\Gamma)\neq\sone(\Omega,\Gamma)$, and
$\sone(\BT,\BG)\neq\sone(\BO,\BG)$.
\end{thm}
\begin{proof}
By a theorem of Brendle \cite{JORG},
assuming CH there exists a set of reals $X$ of size continuum
such that all subsets of $X$ satisfy $\sf P$.
(Recall that $\mathsf{P}=\sone(\BO,\BG)$.)

As $\sf P$ is closed under taking Borel (continuous is enough) images,
we may assume that $X\sbst [0,1]$.
For $Y\sbst [0, 1]$, write $Y+1 = \{y+1 : y\in Y\}$ for the translation
of $Y$ by $1$.
As $|X|=\c$ and only $\c$ many out of the $2^\c$ many subsets of $X$
are Borel, there exists a subset $Y$ of $X$ which
is not $F_\sigma$ neither $G_\delta$.
By a theorem of Galvin and Miller \cite{GM},
for such a subset $Y$ the set $(X\sm Y)\cup (Y+1)$
does not satisfy $\sone(\Omega,\Gamma)$.
Set $A=X\sm Y$ and $B=Y+1$.
Then $A$ and $B$ satisfy $\sone(\BO,\BG)$, and
$A\cup B$ does not satisfy
$\sone(\Omega,\Gamma)=\binom{\Omega}{\Tau}\cap\binom{\Tau}{\Gamma}$.
By Theorem \ref{addT}, $A\cup B$ satisfies $\binom{\Tau}{\Gamma}$
and therefore it does not satisfies $\binom{\Omega}{\Tau}$.
But by Theorem \ref{addS1BTBG}, the set
$A\cup B$ satisfies $\sone(\BT,\BG)$.
\end{proof}

\section{Special elements}

\Subsection{The Cantor set $C$}
Let $C\sbst\R$ be the canonic middle-third Cantor set.

\begin{prop}
Cantor's set $C$ does not satisfy $\sfin(\Gamma,\Tau)$.
\end{prop}
\begin{proof}
Had it satisfied this property, we would have by Theorem \ref{baire-T}
that $C\in\binom{\Tau}{\Gamma}\cap\sfin(\Gamma,\Tau)=\sone(\Gamma,\Gamma)$, contradicting
\cite{coc2}.
\end{proof}
Thus $C$ satisfies $\sfin(\Tau,\Omega)$ and $\ufin(\Gamma,\Tau)$,
and none of the other new properties.

\Subsection{A special Lusin set}
In \cite{huremen2, AddQuad} we construct, using $\cov(\M)=\c$,
special Lusin sets of size $\c$ which satisfy $\sone(\BO,\BO)$.
The meta-structure of the proof is as follows.
At each stage of this construction we define a set $Y_\alpha^*$ which
is a union of less that $\cov(\M)$ many meager sets, and choose
an element $x_\alpha\in G_\alpha\sm Y_\alpha^*$ where $G_\alpha$ is
a basic open subset of $\NZ$.

\begin{thm}
If $\cov(\M)=\c$, then
there exists a Lusin set satisfying $\sone(\BO,\BO)$ but not $\ufin(\Gamma,\Tau)$.
\end{thm}
\begin{proof}
We modify the aforementioned construction
so to make sure that the resulting Lusin set $L$ does not
satisfy the \emp{}. As we do not need to use any group structure,
we will work in $\NN$ rather than $\NZ$.

\begin{lem}
Assume that $A$ is an infinite set of natural numbers, and $f\in\NN$.
Then the sets
\begin{eqnarray*}
       M_{f,A} & = & \{ g\in\NN : [g\le f]\cap A\mbox{ is finite}\}\\
\tilde M_{f,A} & = & \{ g\in\NN : [f < g]\cap A\mbox{ is finite}\}
\end{eqnarray*}
are meager subsets of $\NN$.
\end{lem}
\begin{proof}
For each $k$, the sets
\begin{eqnarray*}
       N_k & = & \{g\in\NN : (\A n>k)\ n\in A\impl f(n)<g(n)\}\\
\tilde N_k & = & \{g\in\NN : (\A n>k)\ n\in A\impl g(n)\le f(n)\}
\end{eqnarray*}
are nowhere dense in $\NN$.
Now, $M_{f,A} = \Union_{k\in\N} N_k$, and
$\tilde M_{f,A} = \Union_{k\in\N} \tilde N_k$.
\end{proof}

Consider an enumeration $\<f_{2\alpha} : \alpha<\c\>$
of $\NN$ which uses only even ordinals.
At stage $\alpha$ for $\alpha$ even, let $Y_\alpha^*$ be
the set defined in \cite{AddQuad}, and let $\tilde Y_\alpha^*$
be the union of $Y_\alpha^*$ and the two meager sets
\begin{eqnarray*}
       M_{f_\alpha,\N} & = & \{ g\in\NN : [g\le f_\alpha]\mbox{ is finite}\}\\
\tilde M_{f_\alpha,\N} & = & \{ g\in\NN : [f_\alpha < g]\mbox{ is finite}\}.
\end{eqnarray*}
Then $\tilde Y_{\alpha^*}$ is a union of less than $\cov(\M)$ many meager sets.
Choose $x_\alpha\in G_\alpha\sm\tilde Y_\alpha^*$.
In step $\alpha+1$ of the construction
let $Y_{\alpha+1}^*$ be defined as in \cite{AddQuad}, and let $\tilde Y_{\alpha+1}^*$
be the union of $Y_{\alpha+1}^*$ with the meager sets
\begin{eqnarray*}
M_{f_\alpha,[f_\alpha < x_\alpha]} & = & \{ g\in\NN : [g\le f_\alpha < x_\alpha]\mbox{ is finite}\}\\
\tilde M_{f_\alpha,[x_\alpha\le f_\alpha]} & = & \{ g\in\NN : [x_\alpha\le f_\alpha<g]\mbox{ is finite}\}.
\end{eqnarray*}
Now choose $x_{\alpha+1}\in G_{\alpha+1}\sm\tilde Y_{\alpha+1}^*$.
Then $x_\alpha$ and $x_{\alpha+1}$ witness that
$f_\alpha$ does not avoid middles in the resulting set $L=\{x_\alpha : \alpha<\c\}$.
Consequently, $L$ does not satisfy $\ufin(\Gamma,\Tau)$.

The proof that $L$ satisfies $\sone(\BO,\BO)$ is as in \cite{AddQuad}.
\end{proof}

\Subsection{Sierpinski sets}
If a Sierpinski set satisfies $\binom{\Tau}{\Gamma}$,
then it satisfies $\sone(\Tau,\Gamma)$ but not
$\sone(\O,\O)$. Such a result would give another
solution to Problem \ref{topoprob}.
However, as towers are null in the usual measure on
${^\N}\{0,1\}$, it is not straightforward to construct
a Sierpinski set which is a member of $\binom{\Tau}{\Gamma}$.
\begin{prob}
Does there exist a Sierpinski set satisfying $\binom{\Tau}{\Gamma}$?
\end{prob}

\Subsection{Unsettled implications}
The paper \cite{CBC} ruled out the possibility
that any selection property for the open case implies any selection property
for the Borel case.
Some implications are ruled out by constructions of \cite{coc2} and \cite{CBC}.
Several other implications are eliminated due to
critical cardinality considerations.

\begin{prob}\label{finalbutone}
Which implications can be added to the diagram
in Figure \ref{survopen} and to the corresponding Borel diagram?
\end{prob}

A summary of all unsettled implications appears in
\cite{SPMprobs}. As a first step towards solving Problem
\ref{finalbutone}, one may try to answer the following.
\begin{prob}
What are the critical cardinalities of the remaining classes?
\footnote{This problem was almost completely
solved in: H.\ Mildenberger, S.\ Shelah, and B.\ Tsaban, \emph{The
combinatorics of $\tau$-covers} (see
\texttt{http://arxiv.org/abs/math.GN/0409068}). There remains exactly
one unsettled critical cardinality in the diagram.}
\end{prob}

\part{Variations on the theme of $\tau$-covers}

\section{$\tau^*$-covers}

The notion of a $\tau^*$-cover is a more flexible
variant of the notion of a $\tau$-cover.

\begin{definition}
A family $Y\sbst\inf$ is \emph{linearly refinable}
if for each $y\in Y$ there exists an infinite subset
$\hat y\sbst y$ such that the family $\hat Y = \{\hat y : y\in Y\}$ is
linearly quasiordered by $\as$.

A cover $\cU$ of $X$ is a \emph{$\tau^*$-cover} of $X$ if it is large,
and $h_\cU[X]$ (where $h_\cU$ is the function defined before Lemma \ref{notions})
is linearly refinable.

For $x\in X$, we will write $x_\cU$ for $h_\cU(x)$, and $\hat x_\cU$ for the
infinite subset of $x_\cU$ such that the sets $\hat x_\cU$ are
linearly quasiordered by $\as$.
\end{definition}

If $\cU$ is a countable $\tau$-cover, then $h_\cU[X]$ is linearly
quasiordered by $\as$ and in particular it is linearly refinable.
Thus every countable $\tau$-cover
is a $\tau^*$-cover. The converse is not
necessarily true.
Let $\Tau^*$ ($\BTstar$) denote the collection of all countable open (Borel)
$\tau^*$-covers of $X$.
Then $\Tau\sbst\Tau^*\sbst\Omega$ and $\BT\sbst\BTstar\sbst\Omega$.

Often problems which are difficult in the case of usual
$\tau$-covers become solvable when shifting to
$\tau^*$-covers. We will give several examples.

\Subsection{Refinements}
One of the major tools in the analysis of selection
principles is to use refinements and de-refinements
of covers. In general, the de-refinement of a
$\tau$-cover is not necessarily a $\tau$-cover.

\begin{lem}\label{tau*refine}
Assume that $\cU\in\Tau^*$ refines
a countable open cover $\cV$ (that is, for each $U\in\cU$ there exists
$V\in\cV$ such that $U\sbst V$). Then $\cV\in\Tau^*$.

The analoguous assertion for countable Borel covers also holds.
\end{lem}
\begin{proof}
Fix a bijective enumeration $\cU=\seq{U_n}$.
Let $\hat x_\cU$, $x\in X$, be as in the definition of $\tau^*$-covers.
For each $n$ let $V_n\in\cV$ be
such that $U_n\sbst V_n$. We claim that $\cW=\{V_n : n\in\N\}\in\Tau^*$.
As $\cW$ is an $\omega$-cover of $X$,
it is infinite;
fix a bijective enumeration $\seq{W_n}$ of $\cW$.
For each $n$ define $S_n = \{k : U_k\sbst W_n\}$,
and $\tilde S_n = S_n\sm\Union_{m<n}S_m$.
For each $x\in X$ define $\hat x_\cW$ by:
$$n\in\hat x_\cW\Iff \tilde S_n\cap\hat x_\cU\neq\emptyset.$$
Then each $\hat x_\cW$ is a subset of $x_\cW$.

Each $\hat x_\cW$ is infinite:
For each
$W_{n_1},\dots,W_{n_k}$
choose
$x_i\nin W_{n_i}$, $i=1,\dots,k$. Then
$\{x,x_1,\dots,x_k\}\not\sbst W_{n_i}$ for all
$i=1,\dots,k$. As $\cU$ is an $\tau^*$-cover of $X$,
there exists $m\in\hat x_\cU$ such that $\{x,x_1,\dots,x_k\}\sbst U_m$.
Consider the (unique) $n$ such that $m\in\tilde S_n$.
Then $U_m\sbst W_n$; therefore
$W_n\nin\{W_{n_1},\dots,W_{n_k}\}$, and
in particular $n\nin\{n_1,\dots,n_k\}$.
As $m\in \tilde S_n\cap \hat x_\cU$, we have that
$n\in \hat x_\cW$.

The sets $\hat x_\cW$ are linearly quasiordered
by $\as$: Assume that $a,b\in X$. We may assume that
$\hat a_\cU\as\hat b_\cU$. As $\lim_n\min\tilde S_n\to\infty$,
we have that $S_n\cap\hat a_\cU\sbst S_n\cap\hat b_\cU$ for all
but finitely many $n$.

This shows that $\cW$ is a $\tau^*$-cover of $X$.
Now, $\cV$ is an extension of $\cW$ by at most countably many elements.
It is easy to see that an extension of a $\tau^*$-cover by
countably many open sets is again a $\tau^*$-cover, see \cite{strongdiags}.
\end{proof}

The first consequence of this important Lemma is that
$\sfin(\fU,\Tau^*)$ ($\sfin(\fU,\BTstar)$) implies
$\ufin(\fU,\Tau^*)$ ($\ufin(\fU,\BTstar)$), that is, the analogue of
Lemma \ref{sfinimplufin} holds.
\begin{cor}
Assume that $\cU=\Union_{n\in\N}\cF_n$, where each $\cF_n$ is finite,
is a $\tau^*$-cover of a space $X$.
Then either $\cup\cF_n = X$ for some $n$, or else
$\cV=\seq{\cup\cF_n}$ is also a $\tau^*$-cover of $X$.
\end{cor}
\begin{proof}
$\cU$ refines $\cV$.
\end{proof}

\Subsection{Equivalences}
All equivalences mentioned in Subsection \ref{equivsec}
hold for $\tau^*$-covers as well.
In particular, the analogue of Theorem \ref{equivalences} holds
(with a similar proof).
\begin{cor}\label{t*g-equiv}
The following equivalences hold:
\be
\i $\sone(\Tau^*,\Gamma)=\sfin(\Tau^*,\Gamma)$;
\i $\sone(\BTstar,\BG)=\sfin(\BTstar,\BG)$.
\ee
\end{cor}
In fact, in the Borel case we get more equivalences in
the case of $\tau^*$-covers than in the case of $\tau$-covers --
see Subsection \ref{t*borelimages}.

\Subsection{Continuous images}
We now solve the problems mentioned in Remarks \ref{clopen-Pwt}
and \ref{ufinGammaTauchar} in the case of $\tau^*$-covers.

\begin{thm}\label{OmTauImage}
The following
properties are equivalent for a set $X$ of reals:
\be
\i $X$ satisfies $\binom{\Omega}{\Tau^*}$;
\i For each continuous image $Y$ of $X$ in $\inf$,
if $Y$ is centered, then $Y$ is linearly refinable.
\ee
\end{thm}
\begin{proof}
$1\Impl 2$: The proof for this is similar to the proof of
$1\Impl 2$ in Theorem \ref{w-tau}.

$2\Impl 1$: Assume that $\cU$ is an $\omega$-cover of $X$.
Replacing each member of $\cU$ with all finite unions of Basic
clopen subsets of it, we may assume that all members of $\cU$
are clopen (to unravel this assumption we will use the fact that
$\Tau^*$ is closed under de-refinements).

Thus, $h_\cU$ is continuous and $Y=h_\cU[X]$ is centered.
Consequently, $Y$ is linearly refinable, that is,
$\cU$ is a $\tau^*$-cover of $X$.
\end{proof}
\begin{rem}\label{BorelOmTauImage}
The analogue assertion (to Theorem \ref{OmTauImage})
for the Borel case,
where open covers are replaced by Borel covers
and continuous image is replaced by Borel image,
also holds and can be proved similarly.
\end{rem}

As in \cite{ShTb}, we will use the notation
$$[f\le h] := \{ n : f(n) \le g(n)\}.$$
Then a subset $Y\sbst\NN$ satisfies the \emp{}
if, and only if, there exists a function $h\in\NN$ such that
the collection
$$\{[f\le h] : f\in Y\}$$
is a subset of $\inf$ and is linearly quasiordered by $\as$.
\begin{definition}
We will say that a subset $Y\sbst\NN$ satisfies the \emph{weak \emp{}}
if there exists a function $h\in\NN$ such that
the collection $\{[f\le h] : f\in Y\}$ is linearly refinable.
\end{definition}

Recall that $\ufin(\Gamma,\Tau^*)=\ufin(\O,\Tau^*)$.
\begin{thm}\label{excludedmiddleopen}
For a zero-dimensional set $X$ of real numbers, the following are equivalent:
\begin{enumerate}
\item{$X$ satisfies $\ufin(\O,\Tau^*)$;}
\item{Every continuous image of $X$ in $\NN$ satisfies the weak \emp{}.}
\end{enumerate}
\end{thm}
\begin{proof}
We make the needed changes in the corresponding proof from \cite{AddQuad}.

$2\Impl 1$: Assume that $\cU_n$, $n\in\N$, are open covers of $X$ which
do not contain finite subcovers.
For each $n$, replacing each member of $\cU_n$
with all of its basic clopen subsets we may assume that all elements of $\cU_n$
are clopen, and thus we may assume further that they are disjoint.
For each $n$ enumerate $\cU_n = \{U^n_m\}_{m\in\N}$.
As we assume that the elements $U^n_m$, $m\in\N$, are disjoint,
we can define a function $\Psi$ from $X$ to $\NN$ by
$$\Psi(x)(n)=m \Iff x\in U^n_m.$$
Then $\Psi$ is continuous. Therefore, $Y=\Psi[X]$
satisfies the weak \emp{}.
Let $h\in\NN$, and for each $f\in Y$, $A_f\sbst [f\le h]$
be such that $\{A_f : f\in Y\}$ is linearly quasiordered by
$\as$.

For each $n$ set
$$\cF_n = \{U^n_k : k\le h(n)\}.$$
We claim that $\cU = \seq{\cup\cF_n}$ is a $\tau^*$-cover
of $X$.
We will use the followig property.
\begin{quote}
\bi
\i[($\star$)]
For each finite subset $F$ of $X$ and each
$n\in\bigcap_{x\in F}A_{\Psi(x)}$,
$F\sbst \cup\cF_n$.
\ei
\end{quote}
Let $\seq{\cup\cF_{k_n}}$ be a bijective enumeration of $\cU$,
and let $f\in\NN$ be such that for each $n$,
$\cup\cF_n = \cup\cF_{k_{f(n)}}$.
For each $x\in X$ set $\hat x_\cU = f[A_{\Psi(x)}]$.
We have the following.

$\hat x_\cU$ is a subset of $x_\cU$:
Assume that $f(n)\in \hat x_\cU$, where
$n\in A_{\Psi(x)}$.
Then $x\in\cup\cF_n=\cup\cF_{k_{f(n)}}$, therefore
$f(n)\in x_\cU$.

$\hat x_\cU$ is infinite:
Assume that $f[A_{\Psi(x)}] = \{f(n_1),\dots,f(n_k)\}$
where $n_1,\dots,n_k\in A_{\Psi(x)}$.
For each $i\le k$ choose $x_i\nin\cup\cF_{k_{f(n_i)}}$,
and set $F=\{x,x_1,\dots,x_k\}$.
Then for all $i\le k$ $F\not\sbst\cup\cF_{k_{f(n_i)}}$.
Choose $n\in\bigcap_{a\in F}A_{\Psi(x)}$.
By property ($\star$), $F\sbst\cup\cF_n=\cup\cF_{k_{f(n)}}$,
therefore $f(n)\nin\{f(n_1),\dots,f(n_k)\}$.
But $n\in A_{\Psi(x)}$, thus $f(n)\in \hat x_\cU$, a contradiction.

As the sets $A_{\Psi(x)}$ are linearly quasiordered by $\as$,
so are the sets $\hat x_\cU = f[A_{\Psi(x)}]$.

$1\Impl 2$:
Since $\Psi$ is continuous, $Y=\Psi[X]$  also satisfies
$\ufin(\O,\Tau^*)$.
Consider the basic open covers $\cU_n = \{U^n_m\}_{m\in\N}$ defined by
$U^n_m = \{f\in Y : f(n) = m\}$. Then there exist finite $\cF_n\sbst\cU_n$, $n\in\N$,
such that either $Y=\cup\cF_n$ for some $n$, or else
$\cV=\{\cup\cF_n : n\in\N\}$ is a $\tau^*$-cover of $Y$.

The first case can be split into two sub-cases: If there exists an infinite
set $A\sbst\N$ such that $Y=\cup\cF_n$, then for each $n\in A$
the set $\{f(n) : f\in Y\}$ is
finite, and we can define
$$h(n) =
\begin{cases}
\max\{f(n) : f\in Y\} & n\in A\\
0 & \mbox{otherwise}
\end{cases}
$$
so that $A\sbst [f\le h]$ for each $f\in Y$, and we
are done.
Otherwise $Y=\cup\cF_n$ for only finitely many $n$, therefore
we may replace each $\cF_n$ satisfying $Y=\cup\cF_n$ with $\cF_n=\emptyset$,
so we are in the second case.

The second case is the interesting one.
$\cV=\{\cup\cF_n : n\in\N\}$ is a $\tau^*$-cover of $Y$ --
fix a bijective enumeration $\seq{\cup\cF_{k_n}}$ of $\cV$ and
witnesses $\hat f_\cV$, $f\in Y$, for that.
Define $h(n) = \max\{m : U^n_m\in\cF_n\}$ for each $n$.
Then the subsets $\{k_n : n\in\hat f_\cV\}$ of $[f\le h]$,
$f\in Y$, are infinite and linearly quasiordered by $\as$.
This shows that $Y$ is linearly refinable.
\end{proof}

\Subsection{Borel images}\label{t*borelimages}
Define the following notion.
\bi
\item[$\sT^*$:]
The set of $X\subseteq\reals$ such that for each linearly refinable Borel image $Y$ of $X$
in $\inf$, $Y$ has a \pi{}.
\ei

By the usual method we get the following.
\begin{lem}\label{T*-char}
$\sT^*=\binom{\BTstar}{\BG}$.
\end{lem}

Clearly $\sT^*$ implies $\sT$.
\begin{lem}
$\non(\sT^*)=\t$.
\end{lem}
\begin{proof}
It is easy to see that $\non(\sT^*)$ is the minimal size of
a linearly refinable family $Y\sbst\inf$ which has no \pi{}.
We will show that $\t\le\non(\sT^*)$. Assume that
$Y\sbst\inf$ is a linearly refinable family of size less than $\t$,
and let $\hat Y$ be a linear refinement of $Y$. As $|\hat Y|\le|Y|<\t$,
$\hat Y$ has a \pi{}, which is in particular a \pi{} of $Y$.
\end{proof}

An application of Lemma \ref{T*-char} and the Cancellation Laws implies the
following.
\begin{thm}\label{Bt*g-char}
$\sone(\BTstar,\BG)=\sT^*\cap\sB$.
\end{thm}

We do not know whether $\sT^*=\sT$.
In particular, we have the following
(recall Theorem \ref{baire-T}).

\begin{prob}
Is it true that every analytic set of reals satisfies $\sT^*$?
Does $\NN\in\sT^*$?
\end{prob}

We do not know whether $\sone(\BG,\BT)=\ufin(\BG,\BT)$ or not.%
\footnote{A negative answer follows from the results in:
H.\ Mildenberger, S.\ Shelah, and B.\ Tsaban, \emph{The
combinatorics of $\tau$-covers} (see
\texttt{http://arxiv.org/abs/math.GN/0409068}).}
This can be contrasted with the following result.

\begin{thm}\label{wX}
For a set $X$ of real numbers, the following are equivalent:
\begin{enumerate}
\item{$X$ satisfies $\sone(\BG,\BTstar)$,}
\item{$X$ satisfies $\sfin(\BG,\BTstar)$,}
\item{$X$ satisfies $\ufin(\BG,\BTstar)$;}
\item{Every Borel image of $X$ in $\NN$ satisfies the weak \emp{}.}
\end{enumerate}
\end{thm}
\begin{proof}
Clearly $1\Impl 2\Impl 3$.

$3\Impl 4$: This can be proved like $1\Impl 2$ in Theorem \ref{excludedmiddleopen}.

$4\Impl 1$:
Assume that $\cU_n=\{U^n_k : k\in\N\}$, $n\in\N$, are Borel $\gamma$-covers of $X$.
We may assume that these covers are pairwise disjoint.
Define a function $\Psi:X\to\NN$ so that for each $x$ and $n$:
$$\Psi(x)(n) = \min\{k : (\A m\ge k)\ x\in U^n_m\}.$$
Then $\Psi$ is a Borel map, and so $Y = \Psi[X]$
satisfies the weak \emp{}.
Let $h\in\NN$ and $A_f\sbst[f\le h]$, $f\in Y$, be witnesses for that.
Set $\cU = \seq{U^n_{h(n)}}$.
For each $x\in X$ set $\hat x_\cU = A_{\Psi(x)}$.
Then $\hat x_\cU$ is infinite and $\hat x_\cU\sbst [\Psi(x)\le h]\sbst x_\cU$
for each $x\in X$, and the sets $\hat x_\cU$ are linearly quasiordered by
$\as$.
\end{proof}

\Subsection{Critical cardinalities}
The argument of Theorem \ref{nontw} implies that
$\non(\sfin(\BTstar,\BO))=\non(\sfin(\Tau^*,\Omega))=\d$.
\begin{thm}
$\non(\sone(\BTstar,\BG))=\non(\sone(\Tau^*,\Gamma))=\t$.
\end{thm}
\begin{proof}
By Theorem \ref{Bt*g-char},
$\non(\sone(\BTstar,\BG)) = \min\{\non(\sT^*),\non(\B)\} = \min\{\t,\b\} = \t$.
On the other hand, $\sone(\Tau^*,\Gamma)$ implies $\sone(\Tau,\Gamma)$, whose
critical cardinality is $\t$.
\end{proof}

Define the following properties.
\bi
\item[$\sX$:]
The set of $X\subseteq\reals$ such that each Borel image of $X$
in $\NN$ satisfies the \emp{}.
\item[$\wX$:]
The set of $X\subseteq\reals$ such that each Borel image of $X$
in $\NN$ satisfies the weak \emp{}.
\ei
Recall that by Theorem \ref{excludedmiddle}, $\ufin(\BG,\BT)=\sX$.
In Theorem \ref{wX} we proved that $\sone(\BG,\BTstar)=\wX$.
We do not know whether $\wX=\sX$.
\begin{prob}
Does $\non(\wX)=\fx$?
\end{prob}

\Subsection{Finite powers}
In \cite{coc2} it is observed that if $\cU$ is an $\omega$-cover of $X$,
then for each $k$ $\cU^k = \{U^k : U\in\cU\}$ is an $\omega$-cover of $X^k$.
Similarly, it is observed in \cite{tau} that
if $\cU$ is a $\tau$-cover of $X$,
then for each $k$ $\cU^k$ is a $\tau$-cover of $X^k$.
We will need the same assertion for $\tau^*$-covers.
\begin{lem}\label{taupow}
Assume that $\cU$ is a $\tau^*$-cover of $X$.
Then for each $k$, $\cU^k$ is a $\tau^*$-cover of $X^k$.
\end{lem}
\begin{proof}
Fix $k$.
Let $\cU=\seq{U_n}$ be an enumeration of $\cU$,
and let $\hat x_\cU\sbst x_\cU$, $x\in X$, witness
that $\cU$ is a $\tau^*$-cover of $X$.
For each $\vec x = (x_0,\dots,x_{k-1})\in X^k$ define
$$\hat{\vec x}_{\cU^k} = \bigcap_{i<k}\hat{(x_i)}_\cU.$$
As the sets $\hat x_\cU$ are infinite and linearly quasiordered
by $\as$, the sets
$\hat{\vec x}_{\cU^k}$ are also infinite and linearly quasiordered
by $\as$. Moreover, for each $n\in \hat{\vec x}_{\cU^k}$
and each $i<k$, $n\in \hat{(x_i)}_\cU$, and therefore
$x_i\in U_n$ for each $i<k$; thus $\vec x\in U_n^k$, as required.
\end{proof}

In \cite{coc2} it is proved that the classes
$\sone(\Omega,\Gamma)$,
$\sone(\Omega,\Omega)$, and
$\sfin(\Omega,\Omega)$ are closed under taking finite
powers, and that none of the remaining classes they considered
has this property.
Actually, their argument for the last assertion
shows that assuming CH, there exist a Lusin set $L$
and a Sierpinski set $S$ such that $L\x L$ and $S\x S$ can be
mapped continuously onto the Baire space $\NN$.
Consequently, we have that none of the classes
$\sone(\Gamma,\Tau)$,
$\sfin(\Gamma,\Tau)$,
$\ufin(\Gamma,\Tau)$,
$\sone(\Tau,\O)$,
and their corresponding Borel versions,
is closed under taking finite powers.
We do not know whether the remaining $7$ classes
which involve $\tau$-covers are closed
under taking finite powers.
\begin{thm}
$\sone(\Omega,\Tau^*)$ and $\sfin(\Omega,\Tau^*)$ are closed
under taking finite powers.
\end{thm}
\begin{proof}
We will prove the assertion for $\sone(\Omega,\Tau^*)$; the proof
for the remaining assertion is similar.
Fix $k$.
In \cite{coc2} it is proved that for each open $\omega$-cover $\cU$ of
$X^k$ there exists an open $\omega$-cover $\cV$ of $X$
such that the $\omega$-cover $\cV^k$ of $X^k$ refines $\cU$.

Assume that $\seq{\cU_n}$ is a sequence of open $\omega$-covers of
$X^k$. For each $n$ choose an open $\omega$-cover $\cV_n$ of $X$
such that $\cV_n^k$ refines $\cU_n$.
Apply $\sone(\Omega,\Tau^*)$ to extract elements $V_n\in\cV_n$, $n\in\N$,
such that $\cW=\seq{V_n}\in\Tau^*$. By Lemma \ref{taupow},
$\cW^k$ is a $\tau^*$-cover of $X^k$.
For each $n$ choose $U_n\in\cU_n$ such that $V_n^k\sbst U_n$.
Then by Lemma \ref{tau*refine}, $\seq{U_n}$ is a $\tau^*$-cover
of $X$.
\end{proof}

\Subsection{Strong properties}
Assume that $\seq{\fU_n}$ is a sequence of collections of covers of a space $X$,
and that $\fV$ is a collection of covers of $X$.
The following selection principle is defined in \cite{strongdiags}.
\begin{itemize}
\item[$\sone(\seq{\fU_n},\fV)$:]
For each sequence $\seq{\cU_n}$ where for each $n$ $\cU_n\in\fU_n$,
there is a sequence
$\seq{U_n}$ such that for each $n$ $U_n\in\cU_n$, and $\seq{U_n}\in\fV$.
\end{itemize}
The notion of \emph{strong $\gamma$-set}, which is due to Galvin and Miller
\cite{GM}, is a particular instance of the new selection principle,
where $\fV=\Gamma$ and for each $n$ $\fU_n=\O_n$, the collection of
open $n$-covers of $X$ (we use here the simple characterization
given in \cite{strongdiags}).
It is well known that the $\gamma$-property $\sone(\Omega,\Gamma)$ does
not imply the strong $\gamma$-property $\sone(\seq{\O_n},\Gamma)$.
It is an open problem whether $\sone(\Omega,\Tau)$ implies
$\sone(\seq{\O_n},\Tau)$.

The following notions are defined in \cite{strongdiags}.
A collection $\fU$ of open covers of a space $X$ is
\emph{finitely thick} if:
\be
\i\label{larger}
If $\cU\in\fU$ and for each $U\in\cU$ $\cF_U$ is a finite
family of open sets such that for each $V\in\cF_U$,
$U\sbst V\neq X$, then $\Union_{U\in\cU}\cF_U\in\fU$.
\i\label{finadd}
If $\cU\in\fU$ and $\cV=\cU\cup\cF$ where $\cF$ is finite
and $X\nin\cF$, then $\cV\in\fU$.
\ee
A collection $\fU$ of open covers of a space $X$ is
\emph{countably thick} if for each $\cU\in\fU$ and
each countable family $\cV$ of open subsets of $X$
such that $X\nin\cV$, $\cU\cup\cV\in\fU$.

Whereas $\Tau$ is in general not finitely thick nor
countably thick, $\Tau^*$ is both finitely and
countably thick \cite{strongdiags}.
In \cite{strongdiags} it is proved that if $\fV$ is
countably thick, then $\sone(\Omega,\fV)=\sone(\seq{\O_n},\fV)$.
Consequently, $\sone(\Omega,\Tau^*)=\sone(\seq{\O_n},\Tau^*)$.

\Subsection{Closing on the Minimal Tower problem}
Clearly $\sone(\Tau^*,\Gamma)$ implies $\sone(\Tau,\Gamma)$,
and $\sone(\BTstar,\BG)$ implies $\sone(\BT,\BG)$.
So we now have new topological lower bounds
on the Minimal Tower problem.
\begin{prob}
\be
\i Does $\sone(\Omega,\Gamma)=\sone(\Tau^*,\Gamma)$?
\i Does $\sone(\BO,\BG)=\sone(\BTstar,\BG)$?
\ee
\end{prob}

We also have a new combinatorial bound.
\begin{definition}
$\p^*$ is the minimal size of a centered family
in $\inf$ which is not linearly refineable.
\end{definition}
\begin{thm}
The critical cardinalities of the properties
$\binom{\BO}{\BTstar}$,
$\sone(\BO,\BTstar)$,
$\sfin(\BO,\BTstar)$,
$\binom{\Omega}{\Tau^*}$,
$\sone(\Omega,\Tau^*)$, and $\sfin(\Omega,\Tau^*)$ are all equal to
is $\p^*$.
\end{thm}
\begin{proof}
By Theorem \ref{OmTauImage} and Remark \ref{OmTauImage},
we have that the assertion holds for $\binom{\Omega}{\Tau^*}$
and $\binom{\BO}{\BTstar}$.

We will use the following result which is analogous to
Theorem \ref{great} and can be proved similarly.
\begin{thm}
The following inclusions hold:
\be
\i $\binom{\Omega}{\Tau^*}\sbst\sfin(\Gamma,\Tau^*)$.
\i $\binom{\BO}{\BTstar}\sbst\sfin(\BG,\BTstar)$.
\ee
\end{thm}
Consequently, $\p^*\le\d$.
As in Theorem \ref{non=kwt}, we get
the remaining assertions follow from this and Corollary \ref{equivalences2},
as $\Tau^*$ is closed under taking countable supersets.
\end{proof}

The following can be proved either directly from the definitions
or from the equivalence
$\binom{\Omega}{\Gamma}=\binom{\Omega}{\Tau^*}\cap\binom{\Tau^*}{\Gamma}$.
\begin{cor}
$\p = \min\{\p^*,\t\}$.
Thus, if $\p<\t$ is consistent, then $\p^*<\t$ is
consistent as well.
\end{cor}

We therefore have the following problem.
\begin{prob}
Does $\p=\p^*$?
\end{prob}

\Subsection{Remaining Borel classes}
We are left with Figure~\ref{survborel} for the Borel case.

\begin{figure}[!h]
\begin{changemargin}{-5cm}{-5cm}
\begin{center}
{\scriptsize
$\xymatrix@C=10pt@R=16pt{
\sone(\BG,\BG)\ar[r]
& \sone(\BG,\BTstar)\ar[rrr]
&
&
& \sone(\BG,\BO)\ar[r]
&\sone(\BG,\B)
\\
&
& \sfin(\BTstar,\BTstar)\ar[ul]\ar[r]
& \sfin(\BTstar,\BO)\ar[ur]
\\
&
& \sfin(\BO,\BTstar)\ar[u]\ar[r]
& \sfin(\BO,\BO)\ar[u]
\\
\sone(\BTstar,\BG)\ar[uuu]\ar[r]
& \sone(\BTstar,\BTstar)\ar[uur]\ar[rrr]
&
&
& \sone(\BTstar,\BO)\ar[uul]\ar[r]
&\sone(\BTstar,\B)\ar[uuu]
\\
\sone(\BO,\BG)\ar[u]\ar[r]
& \sone(\BO,\BTstar)\ar[u]\ar[rrr]
\ar[uur]|!{[u];[ur]}\hole
&
&
&\sone(\BO,\BO)\ar[u]\ar[r]
\ar[uul]|!{[u];[ul]}\hole
& \sone(\B,\B)\ar[u]
\\
}$
}
\caption{\label{survborel}}
\end{center}
\end{changemargin}

\end{figure}

\section{Sequences of compatible $\tau$-covers}\label{LE-covers}

When considering \emph{sequences} of $\tau$-covers,
it may be convenient to have that the linear quasiorderings they define on $X$
agree, in the sense that there exists a Borel linear quasiordering
$\LE$ on $X$ which is contained in all of the induced quasiorderings.
In this case, we say that the $\tau$-covers are \emph{compatible}.
We thus have the following new
selection principle:
\bi
\item[$\sone^\LE(\Tau,\fV)$:]{
For each sequence $\seq{\cU_n}$ of countable open \emph{compatible}
$\tau$-covers of $X$ there is a sequence
$\seq{U_n}$ such that for each $n$ $U_n\in\cU_n$, and $\seq{U_n}\in\fV$.
}
\ei
The selection principle $\sfin^\LE(\Tau,\fV)$ is defined similarly.
Replacing ``open'' by ``Borel'' gives the selection principles
$\sone^\LE(\BT,\fV)$ and $\sfin^\LE(\BT,\fV)$.
The following implications hold:
$$\begin{matrix}
\sone^\LE(\Tau,\fV) & \impl & \sfin^\LE(\Tau,\fV) & \impl & \binom{\Tau}{\fV}\\
\uparrow & & \uparrow\\
\sone(\Tau,\fV) & \impl & \sfin(\Tau,\fV)
\end{matrix}$$
and similarly for the Borel case.

For $\fV=\Gamma$ the new notions coincide with the old ones.
\begin{prop}
The following equivalences hold:
\be
\i $\sone(\Tau,\Gamma)=\sone^\LE(\Tau,\Gamma)
=\sfin(\Tau,\Gamma)=\sfin^\LE(\Tau,\Gamma)$;
\i $\sone(\BT,\BG)=\sone^\LE(\BT,\BG)=
\sfin(\BT,\BG)=\sfin^\LE(\BT,\BG)$.
\ee
\end{prop}
\begin{proof}
We will prove (1); (2) is similar.
By Theorem \ref{equivalences}, we have the following implications
$$\begin{matrix}
\sone^\LE(\Tau,\Gamma) & \impl & \sfin^\LE(\Tau,\Gamma) & \impl &
\binom{\Tau}{\Gamma}\cap\sfin(\Gamma,\Gamma)=\sone(\Tau,\Gamma)\\
\uparrow & &\uparrow\\
\sone(\Tau,\Gamma) & \impl & \sfin(\Tau,\Gamma)
\end{matrix}$$
\end{proof}

\Subsection{The class $\sfin^\LE(\BT,\BT)$}

\begin{definition}
A \emph{$\tau$-cover of $\<X,\LE\>$} is a $\tau$-cover of $X$ such that
the induced quasiordering contains $\LE$.
\end{definition}

\begin{lem}\label{square}
Let $\<X,\LE\>$ be a linearly quasiordered set of reals, and
assume that
every Borel image of $\LE$ in $\NN$ is bounded (with respect to $\le^*$).
Assume that $\cU_n=\{U^n_k : k\in\N\}$ are Borel $\tau$-covers of $\<X,\LE\>$.
Then there exist finite subsets $\cF_n$ of $\cU_n$, $n\in\N$,
such that $\Union_{n\in\N}\cF_n$ is a $\tau$-cover of $\<X,\LE\>$.
\end{lem}
\begin{proof}
Fix a linear quasiordering $\LE$ of $X$, and
assume that $\cU_n=\{U^n_k : k\in\N\}$ are Borel $\tau$-covers of $\<X,\LE\>$.
Define a Borel function $\Psi$ from $\LE$ to $\NN$ by:
$$\Psi(x,y)(n) = \min\{k : \(\forall m\geq k\)\ x\in U^n_m\impl y\in U^n_m\}.$$
$\Psi[\LE]$ is bounded, say by $g$.
Now define a Borel function $\Phi$ from $X$ to $\NN$ by:
$$\Phi(x)(n) = \min\{k : g(n)\le k\mathrm{\ and \ } x\in U^n_k\}$$
Note that $\Phi[X]$ is a Borel image of $\LE$ in $\NN$,
thus it is bounded, say by $f$.
It follows that the sequence $\seq{U^n_{g(n)},\dots,U^n_{f(n)} : g(n)\le f(n)}$
is large, and is a $\tau$-cover of $X$.
\end{proof}

According to \cite{CBC},
the property that every Borel image is bounded is equivalent to $\sone(\BG,\BG)$.

\begin{lem}\label{property}
Let ${\mathcal P}$ be a collection of spaces which is
closed under taking Borel subsets, continuous images (or isometries),
and finite unions.
Then for each set $X$ of real numbers, the following are equivalent:
\be
\i Each Borel linear quasiordering $\LE$ of $X$ satisfies $\mathcal P$,
\i $X^2$ satisfies $\mathcal P$.
\i There exists a Borel linear quasiordering $\LE$ of $X$ satisfying $\mathcal P$,
\ee
\end{lem}
\begin{proof}
$1\Impl 2\Impl 3$: The set $\LE=X^2$ is a linear quasiordering of $X$.

$3\Impl 2$: If $\LE$ satisfies $\mathcal P$, then so does
its continuous image $\GE=\{(y,x) : x\LE y\}$.
Thus, $X^2=\GE\cup\LE$ satisfies $\mathcal P$.

$2\Impl 1$: $\mathcal P$ is closed under taking Borel subsets.
\end{proof}

Thus, lemma \ref{square} can be restated as follows.

\begin{thm}
If $X^2$ satisfies $\sone(\BG,\BG)$, then
$X$ satisfies $\sfin^\LE(\BT,\BT)$.
\end{thm}
\begin{proof}
The property $\sone(\BG,\BG)$ satisfies the assumptions of Lemma \ref{property}.
\end{proof}

\begin{prob}\label{prob1}
Assume that $X$ satisfies $\sfin^\LE(\BT,\BT)$.
Is it true that $X^2$ satisfies $\sone(\BG,\BG)$?
\end{prob}

Since $\b = \cov({\mathcal N})=\cof({\mathcal N})$ (in particular, the Continuum Hypothesis)
implies that $\sone(\BG,\BG)$ is not closed under taking squares \cite{CBC}, a
positive answer to Problem \ref{prob1} would imply that
the property $\sone(\BG,\BG)$ does not imply $\sfin^\LE(\BT,\BT)$.

\bigskip

{\small \textbf{Acknowledgments.}
We wish to thank Martin Goldstern for reading this paper
and making several important
comments.
We also thank
Saharon Shelah for the opportunity to introduce this work
at his research seminar, and for the fruitful
cooperation which followed \cite{ShTb}.
Finally, we would like to thank the organizers of the
\emph{Lecce Workshop on Coverings, Selections and Games in Topology}
(June 2002)---Cosimo Guido, Ljubisa D.R.\ Ko\v{c}inac, and Marion Scheepers---for
inviting us to give a lecture on this paper in their workshop and for the kind
hospitality during the workshop.
}



\begin{thebibliography}{00}


\bibitem{huremen2}
T.\ Bartoszy\'nski, S.\ Shelah, and B.\ Tsaban,
\emph{Additivity properties of topological diagonalizations},
The Journal of Symbolic Logic \textbf{68} (2003),
1254--1260.

\bibitem{HBK}
A.\ R.\ Blass,
\emph{Combinatorial cardinal characteristics of the continuum},
in: \textbf{Handbook of Set Theory} (M.\ Foreman, A.\ Kanamori, and M.\ Magidor, eds.),
Kluwer Academic Publishers, Dordrecht,
to appear.

\bibitem{JORG}
J.\ Brendle,
\emph{Generic constructions of small sets of reals},
Topology Appl.\ \textbf{71} (1996), 125--147.

\bibitem{vD}
E.\ K.\ van~Douwen,
\emph{The  integers  and  topology},
in: \textbf{Handbook of Set Theoretic Topology}
(K.\ Kunen  and  J.\ Vaughan, Eds.),
North-Holland, Amsterdam, 1984, 111--167.

\bibitem{GM}
F.\ Galvin and A.W.\ Miller,
\emph{$\gamma$-sets and other singular sets of real numbers},
Topology Appl.\ \textbf{17} (1984), 145--155.

\bibitem{gerlitsnagy}
J.\ Gerlits and Zs.\ Nagy,
\emph{Some properties of C(X), I},
Topology and its Applications \textbf{14} (1982),
151--161.



\bibitem{coc2}
W.\ Just, A.W.\ Miller, M.\ Scheepers, and P.J.\ Szeptycki,
\emph{The combinatorics of open covers II},
Topology and its Applications \textbf{73} (1996), 241--266.

\bibitem{pawlikowskireclaw}
J.\ Pawlikowski and I.\ Rec{\l}aw,
\emph{Parametrized Cicho\'n's diagram and small sets},
Fundamenta    Mathematicae \textbf{147} (1995),
135--155.

\bibitem{reclaw}
I.\ Rec\l{}aw,
\emph{Every Luzin set is undetermined in the point-open game},
Fund.\ Math.\ \textbf{144} (1994),
43--54.


\bibitem{ROTH2}
F.\ Rothberger,
\emph{On some problems of Hausdorff and of Sierpi\'nski},
Fund.\ Math.\ \textbf{35} (1948), 29--46.

\bibitem{coc1}
M.\ Scheepers,
\emph{Combinatorics of open covers I: Ramsey theory},
Topology and its Applications \textbf{69} (1996),
31--62.

\bibitem{wqn}
M.\ Scheepers,
\emph{Sequential convergence in ${\sf C}_p(X)$ and a covering property},
East-West Journal of Mathematics \textbf{1} (1999),
207--214.

\bibitem{CBC}
M.\ Scheepers and B.\ Tsaban,
\emph{The combinatorics of Borel covers},
Topology and its Applications \textbf{121} (2002),
357--382.
\\ \arx{math.GN/0302322}

\bibitem{ShTb}
S.\ Shelah and B.\ Tsaban,
\emph{Critical cardinalities and additivity properties of combinatorial notions of smallness},
Journal of Applied Analysis \textbf{9} (2003),
149--162.
\\ \arx{math.LO/0304019}

\bibitem{tau}
B.\ Tsaban,
\emph{A topological interpretation of $\mathfrak{t}$},
Real Analysis Exchange \textbf{25} (1999/2000),
391--404.
\\ \arx{math.LO/9705209}

\bibitem{huremen}
B.\ Tsaban,
\emph{A diagonalization property between Hurewicz and Menger},
Real Analysis Exchange \textbf{27} (2001/2002),
757--763.
\\ \arx{math.GN/0106085}

\bibitem{SPMprobs}
B.\ Tsaban,
\emph{Selection principles in Mathematics: A milestone of open problems},
Note di Matematica \textbf{22} (2003),
179--208.
\\ \arx{math.GN/0312182}

\bibitem{strongdiags}
B.\ Tsaban,
\emph{Strong $\gamma$-sets and other singular spaces},
Topology and its Applications \textbf{153} (2005),
620--639.

\bibitem{AddQuad}
B.\ Tsaban,
\emph{Additivity numbers of covering properties},
in: \textbf{Selection Principles and Covering Properties in Topology} (L.\ D.R.\ Ko\v{c}inac, ed.),
Quaderni di Matematica,
to appear.
\\ \arx{math.GN/0604451}

\end{thebibliography}
\end{document}